\documentclass{amsart}

\usepackage{amsmath, amssymb, amsthm, amsfonts}
\usepackage{MnSymbol}
\usepackage{multicol}
\usepackage{tikz-cd, xcolor}
\usepackage{hyperref, caption}
\usepackage{ enumerate, adjustbox}
\usepackage{todonotes}
\hypersetup{
	colorlinks,
	citecolor=green,
	filecolor=black,
	linkcolor=blue,
	urlcolor=black}
\usepackage[utf8]{inputenc}
\usepackage[english]{babel}
\usepackage{csquotes}
\usepackage[backend=biber,
style=alphabetic,
sorting=nyt,maxbibnames=9, maxcitenames=9]{biblatex}
\newcommand{\Z}{\mathbb{Z}}

\newcommand{\R}{\mathbb{R}}

\DeclareMathOperator{\colim}{colim}

\newcommand{\mM}{\underline{M}}
\newcommand{\mR}{\underline{R}}

\newcommand{\mN}{\underline{N}}
\newcommand{\mpi}{\underline{\pi}}
\newcommand{\mhr}{\underline{\hr}}
\newcommand{\mhh}{\underline{\hh}}

\newcommand{\mF}{\underline{\mathbb{F}}}
\newcommand{\smsh}{\wedge}
\newcommand{\boxover}[1]{\ensuremath{\mathop{\Box}_{#1}}}

\newcommand{\bok}{B{\"o}kstedt}
\newcommand{\bord}{MU_{\R}}
\newcommand{\zetanorm}{N^{D_{2m}}_{\zeta D_2 \zeta^{-1}}c_\zeta}
\newcommand{\simpotwo}{O(2)_{\bullet}}
\DeclareMathOperator{\thh}{THH}
\DeclareMathOperator{\thr}{THR}
\DeclareMathOperator{\hh}{HH}
\DeclareMathOperator{\hr}{HR}

\theoremstyle{plain} 
\newtheorem*{thm*}{Theorem}
\newtheorem{thm}{Theorem}[section]
\newtheorem{cor}[thm]{Corollary}
\newtheorem{lemma}[thm]{Lemma}
\newtheorem{prop}[thm]{Proposition}

\theoremstyle{definition}
\newtheorem{defn}[thm]{Definition}
\newtheorem*{nota}{Convention}

\theoremstyle{remark}
\newtheorem{rem}[thm]{Remark}
\newtheorem{ex}[thm]{Example}

\numberwithin{equation}{section}
\numberwithin{figure}{section}

\title{Computational tools for Real topological Hochschild homology}
\author{Chloe Lewis}
\address{Department of Mathematics, University of Wisconsin-Eau Claire, Eau Claire, WI 54701}
\email{lewischl@uwec.edu}
\date{August 2024}

\begin{document}

\begin{abstract}
In this paper, we construct a Real equivariant version of the \bok\ spectral sequence which takes inputs in the theory of Real Hochschild homology developed by Angelini-Knoll, Gerhardt, and Hill and converges to the equivariant homology of Real topological Hochschild homology, $\thr$. We also show that when $A$ is a commutative $C_2$-ring spectrum, $\thr(A)$ has the structure of an $A$-Hopf algebroid in the $C_2$-equivariant stable homotopy category.

\end{abstract}

\maketitle


\section{Introduction}

A very successful approach in the modern study of algebraic $K$-theory is that of trace methods. Historically, the algebraic $K$-theory groups of a ring $A$ have proven to be quite difficult to compute; the trace methods approach uses more computationally accessible invariants of rings (and their topological analogues) as approximations of $K$-theory. Hochschild homology, denoted by $\hh$, is one such approximation via the Dennis trace map, 
\[
K_n(A) \rightarrow \hh_n(A).
\]
Hochschild homology is a well-studied invariant in homological algebra and therefore computations are relatively approachable. However, the Dennis trace does not provide a particularly good approximation of $K$-theory. This is remedied by using a topological version of Hochschild homology which also receives a trace map from $K$-theory. Topological Hochschild homology, denoted by $\thh$, is also a necessary building block in topological cyclic homology, which has proven to be a very good approximation of $K$-theory and plays an important role in many modern $K$-theory computations. 

Work of \bok\ \cite{bokstedtthhZandZp} gives a spectral sequence which takes input data from classical Hochschild homology groups and produces information about the homology of $\thh$. For a ring spectrum $R$ and a field $k$, this spectral sequence takes the form 
\[
    E^2_{*,*}= \hh_*(H_*(R);k) \Rightarrow H_{*}(\thh(R);k).
\]
The \bok\ spectral sequence has proven quite useful in $\thh$ calculations since it leverages the well-understood theory of classical Hochschild homology to generate information about $\thh$. 

The story of trace methods for algebraic $K$-theory can be generalized to an equivariant setting by considering invariants of rings and ring spectra endowed with a group action. One such equivariant generalization is Real algebraic $K$-theory, denoted by $KR$, for rings with the $C_2$-action of involution. Real $K$-theory, which also generalizes Hermitian $K$-theory \cite{karoubi}, was originally defined by Hesselholt-Madsen in \cite{realktheory}. Concurrently, the authors developed an equivariant version of $\thh$ called Real topological Hochschild homology, denoted by $\thr$. In the Real equivariant setting, $\thr$ receives a map from Real algebraic $K$-theory, and an analogous story of approximating $KR$ via Real topological cyclic homology and $\thr$ unfolds (see \cite{dottothesis,dottosolo,hermitian,hogenhaven}).   

A definition of $\thr$ via a simplicial bar construction was given by Dotto, Moi, Patchkoria, and Reeh in \cite{dotto}. This construction lends itself particularly well to further computational work. Recently, Angelini-Knoll, Gerhardt, and Hill \cite{realhh} showed that $\thr(A)$ can be modeled as the equivariant norm $N_{C_2}^{O(2)}(A)$ and constructed an algebraic analogue to $\thr$ which they call Real Hochschild homology. In this paper, we extend this definition of Real Hochschild homology to the graded setting.  

The main result of this paper is the construction of a Real \bok\ spectral sequence which computes the dihedral $D_{2m}$-equivariant homology of $\thr$ using input data from Real $D_{2m}$-Hochschild homology. This result is restated as Theorem \ref{thm: real bok ss} in the text.  

\begin{thm} 
     Let $A$ be a ring spectrum with anti-involution and let $E$ be a commutative $D_{2m}$-ring spectrum. If $\underline{E}_{\ostar}(N^{D_{2m}}_{C_2} A)$ and $\underline{E}_{\ostar}(N^{D_{2m}}_e \iota^*_e A)$ are both flat as modules over $\underline{E}_{\ostar}$ and if $A$ has free $(\iota_{C_2}^*E)$- and $\iota_{e}^*E$-homology then there is a Real \bok\ spectral sequence of the form
    \[
    E^2_{*, \ostar} = \mhr^{\underline{E}_{\ostar},D_{2m}}_*(\underline{(\iota_{C_2}^*E)}_{\ostar}(A)) \Rightarrow \underline{E}_{\ostar}(\iota_{D_{2m}}^* \thr(A)).
    \]
\end{thm}
\noindent Here, $E_{\ostar}$ denotes an $\underline{RO}(D_{2m})$-graded homology theory. We review $\underline{RO}$-gradings in Section \ref{section: graded inputs for hh}. Of particular interest is the case of $E=H\mF_2$ where we obtain the following corollary:

\begin{cor}
    Let $A$ be a ring spectrum with anti-involution and such that $\underline{H}^{C_2}_{\ostar}(A; \mF_2)$ and $\underline{H}^{C_2}_{\ostar}(N^{C_2}_e \iota_e^* A; \mF_2)$ are flat as modules over ${H\mF_2}_{\ostar}$. Then there is a Real \bok\ spectral sequence
    \[
    E^2_{*, \ostar}=\mhr^{C_2}_*(\underline{H}_{\ostar}^{C_2} (A; \mF_2)) \Rightarrow \underline{H}_{\ostar}^{C_2}(\thr(A); \mF_2).
    \]
\end{cor}

The techniques used to construct the Real \bok\ spectral sequence may be extended to a different flavor of equivariant topological Hochschild homology. For $G$ a finite subgroup of $S^1$, the $G$-twisted topological Hochschild homology of a ring spectrum $R$, denoted by $\thh_G(R)$, is an $S^1$-spectrum which incorporates an action of the cyclic group into the constructions. In particular, the norm $N^{S^1}_G R$ is a model for $\thh_G(R)$. In this twisted setting there is an analogous algebraic theory of Hochschild homology for Green functors. Work of Adamyk-Gerhardt-Hess-Klang-Kong \cite{twistedthh} constructs a $G$-twisted \bok\ spectral sequence which computes the $G$-equivariant homology of $\thh_G$. In this paper, we construct a spectral sequence which computes the $G$-equivariant homology of $H$-twisted $\thh$, for a subgroup $H$ of $G$. This result is restated as Theorem \ref{thm: relative twisted bok ss} in Section \ref{section: higher twisted bok ss}. 

\begin{thm}
     Let $H \subseteq G$ be finite subgroups of $S^1$ and let $g=e^{2\pi i/|G|}$ be a generator of $G$. Let $R$ be an $H$-ring spectrum and $E$ a commutative $G$-ring spectrum. Assume that $g$ acts trivially on $E$ and that $\underline{E}_{\ostar}(N^G_H R)$ is flat as a module over $\underline{E}_{\ostar}$. If $R$ has $(\iota_H^*E)$-free homology, then there is a relative twisted \bok\ spectral sequence
      \[
      E_{*,\ostar}^2={\mhh_H^{\underline{E}_{\ostar},G}}(\underline{(\iota_H^*E)}_{\ostar}(R))_* \Rightarrow \underline{E}_{*+\ostar}(\iota^*_G\thh_H(R)).
      \]
\end{thm}
Taking $G=H$ in this theorem recovers the spectral sequence of \cite{twistedthh}. In the case of $H=e$, this result also gives a new spectral sequence converging to the $G$-equivariant homology of ordinary $\thh$. 

One way to gain computational traction in calculations involving the classical \bok\ spectral sequence is to utilize the algebraic structures present in $\thh$ and in the spectral sequence itself. Work of Angeltveit-Rognes \cite{angeltveit2005hopf} shows that, when the input is commutative and under appropriate flatness conditions, the Hopf algebra structure in $\thh$ lifts to the \bok\ spectral sequence. In this paper, we prove that when the input is commutative, $\thr$ has the structure of a Hopf algebroid in the $C_2$-equivariant stable homotopy category. This result is restated as Theorem \ref{thm: thr is a hopf algebroid} in Section \ref{chapter: algebraic structures}.  
\begin{thm}
    Let $A$ be a commutative $C_2$-ring spectrum. The Real topological Hochschild homology of $A$ is a Hopf algebroid in the $C_2$-equivariant stable homotopy category. 
\end{thm}
In future work, we show that this Hopf algebroid structure lifts to the Real \bok\ spectral sequence. 

\subsection{Organization}
In Section \ref{section: equivar thh} we review Real Hochschild theories in the topological and algebraic settings. We also define a new notion of graded Real Hochschild homology in \ref{section: graded inputs for hh}. In Section \ref{section: construction} we recall Hill's notion of free homology and use it to construct the Real \bok\ spectral sequence. In Section \ref{section: twisted thh} we construct an analogous spectral sequence for $G$-twisted topological Hochschild homology. Finally, in Section \ref{chapter: algebraic structures} we prove that $\thr(A)$ has the structure of a Hopf algebroid when $A$ is a commutative $C_2$-ring spectrum.  

\subsection{Conventions}
In the remainder of this paper we will use the notation $D_2$ to denote the cyclic group of two elements. A star $*$ refers to an integer grading, the star $\star$ is used to denote an $RO(G)$-grading, and the symbol $\ostar$ is used for an $\underline{RO}(G)$-grading. Simplicial gradings are denoted by $\bullet$.

\subsection{Acknowledgements}
This work first appeared in the author's Ph.D. thesis. She would like to thank her thesis advisor Teena Gerhardt for her helpful guidance and excellent mentorship on this project. The author would also like to thank Mike Hill and Danika Van Niel for several helpful conversations, particularly those pertaining to the Hopf algebroid structure.

This work was supported by NSF grants DMS-2104233 and DMS-1810575 as well as the Graduate School and Mathematics Department at Michigan State University. 


\section{Real (topological) Hochschild homology}\label{section: equivar thh}

We begin by recalling the Real equivariant variations of Hochschild homology in both the algebraic and topological settings. Classically, Hochschild homology and topological Hochschild homology can be defined via simplicial constructions and we begin this section by reviewing the appropriate notion of simplicial objects in the Real equivariant setting. 


\subsection{Equivariant simplicial objects} \label{subsection: equivariant simplicial objects}

To set notation, let $D_{2m}$ denote the dihedral group generated by two elements,
\[
D_{2m}= \langle \omega, t \mid\ \omega^2=t^m=1, \omega t \omega=t^{-1} \rangle.
\]
Note that when $m=1$, this is the cyclic group on two elements typically denoted by $C_2$. In anticipation of our discussion about dihedral-equivariant objects in algebra and topology, we elect to use the notation $D_2$ for this group instead. 

\begin{defn} \label{defn: dihedral object}
A \emph{dihedral object} $L_{\bullet}$ in a category $\mathcal{C}$ is a simplicial object in $\mathcal{C}$ together with a $D_{2(n+1)}$-action on $L_n$ specified by the action of the generators:
\[
    t_n: L_n \rightarrow L_n\ \text{ and }\ \omega_n: L_n \rightarrow L_n
\]
such that the following relations hold:
\begin{multicols}{2}
\begin{enumerate}
    \item $\omega_n \circ t_n = t_n^{-1} \circ \omega_n$
    \item $d_0 \circ t_{n+1} = d_n$
    \item $d_i \circ t_n  = t_{n-1} \circ d_{i-1}$ if $1 \leq i \leq n$
    \item $s_0 \circ t_n = t^2_{n+1} \circ s_n$
    \item $s_i \circ t_n = t_{n+1} \circ s_{i-1}$ if $1 \leq i \leq n$
    \item $d_i \circ \omega_n = \omega_{n-1}\circ d_{n-i}$ if $0 \leq i \leq n$
    \item $s_i \circ \omega_n = \omega_{n+1}\circ s_{n-i}$ if $0 \leq i \leq n$.
\end{enumerate}
\end{multicols}
\end{defn}

\begin{defn} \label{defn: real simplicial object}
A \emph{Real simplicial object} $K_{\bullet}$ is a simplicial object together with maps $\omega_n: K_n \rightarrow K_n$ for each $n \geq 0$ which square to the identity and obey relations 6 and 7 in Definition \ref{defn: dihedral object}.  
\end{defn}

\begin{rem}
    By Theorem 5.3 of \cite{fiedorowicz1991crossed}, the geometric realization of a dihedral object has an action of the orthogonal group $O(2)$ and the geometric realization of a Real simplicial object has a $D_2$-action. 
\end{rem}

At times, it is useful to work with a subdivided simplicial object which supports a $D_2$-action. The appropriate subdivision in this case is attributed to Segal \cite{segal1973configuration} and Quillen.  

\begin{defn} \label{defn: Segal subdivision}
    The \emph{Segal-Quillen subdivision} of a simplicial object $X_{\bullet}$, denoted $\text{sq} X_{\bullet}$, is the simplicial object with $k$-simplices
    \[
    \text{sq}  X_k = X_{2k+1}.
    \]
    Let $d_i$ and $s_i$ denote the face and degeneracy maps of $X_{\bullet}$. The face and degeneracy maps $\tilde{d_i}$ and $\tilde{s_i}$ of the subdivision sq$X_{\bullet}$ are given by 
    \begin{align*}
    \tilde{d_i} &= d_i \circ  d_{2k+1-i}\\
    \tilde{s_i} &= s_{2k-i} \circ s_i.
    \end{align*}
    The geometric realization of a simplicial object and its Segal-Quillen subdivision are homeomorphic by Lemma 2.4 of \cite{spalinski2000homotopy}. Furthermore, if $X_{\bullet}$ is a dihedral or Real simplicial set then the homeomorphism $|X_{\bullet}| \cong |\text{sq} X_{\bullet}|$ is a $D_2$-homeomorphism (see \cite{realhh}, Section 2.1). 
\end{defn}


\subsection{Real topological Hochschild homology}

A definition of Real topological Hochschild homology ($\thr$), was first given by Hesselholt-Madsen in \cite{realktheory} in the style of \bok's original construction of $\thh$ \cite{bokthh}. Dotto-Moi-Patchkoria-Reeh in \cite{dotto} subsequently gave a construction of $\thr$ using a dihedral bar construction analogous to the definition of $\thh$ via the cyclic bar construction.  We recall this definition of $\thr$ beginning with a formal description of its input. 

\begin{defn} \label{defn: ring spectrum with anti-involution}
    A \emph{ring spectrum with anti-involution} is a pair $(A, \omega)$ consisting of a ring spectrum $A$ and a map $\omega: A^{op} \rightarrow A$ such that $\omega^2=id$. Here, $A^{op}$ denotes the ring spectrum $A$ with the reversed multiplicative structure where two factors are permuted before being multiplied. 
\end{defn}

\begin{ex}
    Let $A$ be a commutative $D_2$-ring spectrum. Since $A$ is commutative, $A^{op}=A$ and the $D_2$-action on $A$ defines an anti-involution. 
\end{ex}

Classical $K$-theory takes inputs in rings (or ring spectra) hence the appropriate input to Real algebraic $K$-theory is a ring (or ring spectrum) with anti-involution. A ring $A$ with an anti-involution is a ring equipped with a map $\omega: A^{op} \rightarrow A$ which squares to the identity. The motivating example is a group ring $\Z[G]$ with a $D_2$-action given by inversion. For $g, g' \in G$, applying the action to a product yields $(gg')^{-1}=g'^{-1}g^{-1}$. The reversed multiplicative structure here exemplifies why this map lands in the opposite ring. 

In the bar construction which defines $\thr(A)$ we will need an appropriate choice of an $A$-bimodule to play the role of coefficients. 

Let $(A,\omega)$ be a ring spectrum with anti-involution and $M$ an $A$-bimodule with left action map $\psi_L$ and right action map $\psi_R$. We may define an $A$-bimodule structure on $M^{op}$ via the composite maps
\[
    A \smsh M \xrightarrow{\tau} M \smsh A \xrightarrow{id \smsh \omega} M \smsh A \xrightarrow{\psi_R} M
\]
\[
    M \smsh A \xrightarrow{\tau} A \smsh M \xrightarrow{\omega \smsh id} A \smsh M \xrightarrow{\psi_L} M,
\]
 where $\tau$ denotes the swap map that permutes the two factors. 
An \emph{$(A, \omega)$-bimodule} is a pair $(M, \sigma)$ which consists of an $A$-bimodule $M$ and a map of $A$-bimodules $\sigma: M^{op} \rightarrow M$ such that $\sigma^2=id$.

We now recall the definition of Real topological Hochschild homology via a bar construction due to \cite{dotto}. 

\begin{defn} \label{defn: dihedral bar construction}

Let $(A, \omega)$ be a ring spectrum with anti-involution and $(M, \sigma)$ an $(A, \omega)$-bimodule. The \emph{dihedral bar construction of $(A, \omega)$ with coefficients in $(M, \sigma)$} is a Real simplicial spectrum denoted by $B_{\bullet}^{di}(A;M)$. This spectrum has $k$-simplices
\[
B_k^{di} (A;M) = M \smsh A^{\smsh k}.
\]

The simplicial structure maps in this spectrum are the same as those in the classical cyclic bar construction which defines $\thh$.

Additionally, the dihedral bar construction has a level-wise involution $W$. At level $k$, let $\mathbf{k}$ be the $D_2$-set of integers $k=\{1, ..., k\}$ with a $D_2$-action of $\gamma(i)=k+1-i$ for $\gamma$ a generator of $D_2$. The involution is given by 
\[
W: M \smsh A^{\smsh \mathbf{k}} \xrightarrow{id \smsh A_{\gamma(1)} \smsh ... \smsh A_{\gamma(k)}} M \smsh A^{\smsh \mathbf{k}} \xrightarrow{\sigma \smsh \omega^{\smsh \mathbf{k}}} M \smsh A^{\smsh \mathbf{k}}. 
\]
\end{defn}

\begin{defn} \label{defn: thr}
The \emph{Real topological Hochschild homology} of a ring spectrum with anti-involution $(A, \omega)$  with coefficients in the bimodule $(M,\sigma)$ is the $D_2$-spectrum given by the geometric realization of the dihedral bar construction,
\[
    \thr(A;M) := |B_{\bullet}^{di}(A;M)|.
\]
If taking coefficients in $A$, we omit them from the notation and write simply $\thh(A)$. 
\end{defn}

\begin{rem} \label{remark: thr is o2 spectrum}
  Recall that the cyclic operators at each level of the classical cyclic bar construction realized to give $\thh$ an $S^1$-action. The dihedral bar construction inherits much of the same structure as the cyclic bar construction, including these cyclic operators. The level-wise involution produces a new $D_2$-action giving $B^{di}_{n}(A)$ a $C_n \rtimes D_2=D_{2n}$-action. These actions assemble to an action of $O(2)$ on the geometric realization, thus $\thr(A)$ is an $O(2)$-spectrum.  
\end{rem}

\subsection{Real Hochschild homology}
We now consider the analogous algebraic theory for inputs with an involution. A Hochschild homology theory for rings and algebras equipped with an anti-involution called dihedral homology is described in Section 5.2 of \cite{loday}, however this theory is does not fully capture the Real equivariant structure present. The failure is visible when one attempts to construct a Real linearization map from $\thr$. Recall that classically, a connective link between $\thh$ and $\hh$ is the linearization map 
\[\pi_n(\thh(R)) \rightarrow \hh_n(\pi_0(R)).\]

Given a ring spectrum with anti-involution $A$, a Real linearization map should land in the Real Hochschild homology of $\pi_0(A)$. This is not simply a ring with involution; since $A$ is a $D_2$-spectrum, its homotopy forms the $D_2$-Mackey functor $\mpi_n^{D_2}(A)$. Thus, a true algebraic analogue of $\thr$ must be constructed using the language of equivariant algebra and take inputs in a $D_2$-Mackey functor that encodes the action of involution.

Such a construction is given by Angelini-Knoll, Gerhardt, and Hill in \cite{realhh}. We begin with a definition of the appropriate input for this Real Hochschild theory, a particular kind of $D_2$-Mackey functor called a discrete $E_{\sigma}$-ring. 

\begin{defn} \label{defn: discrete e sigma ring}
   A \emph{discrete $E_{\sigma}$-ring} consists of the following data: 
   \begin{enumerate}
       \item A $D_2$-Mackey functor $\mM$ such that there is an associative product on $\mM(D_2/e)$  for which the Weyl action is an anti-homomorphism.
       \item An $N^{D_2}_e \iota_e^* \mM$-bimodule structure on $\mM$ with right action 
       \[\psi_L: N^{D_2}_e \iota_e^* \mM\ \square\ \mM \rightarrow \mM\] 
       and left action 
       \[\psi_R: \mM\ \square\ N^{D_2}_e \iota_e^* \mM \rightarrow \mM. \] 
       We further require that $\psi$ restricts to the usual module action over the enveloping algebra on $\mM(D_2/e)$.
       \item A unit element $1 \in \mM(D_2/D_2)$ such that $res(1)=1 \in \mM(D_2/e)$. 
   \end{enumerate}
\end{defn}

\begin{rem}
    The data of a discrete $E_{\sigma}$-ring is essentially that of a Hermitian Mackey functor plus the unitality condition in (3) above. The definition of Hermitian Mackey functors is due to Dotto-Ogle; see \cite{hermitian}. 
\end{rem}

\begin{ex}
    Let $A$ be a ring spectrum with anti-involution. Then $\mpi_0^{D_2}(A)$ has the structure of a discrete $E_{\sigma}$-ring. 
\end{ex}

The data of a discrete $E_{\sigma}$-ring $\mM$ includes a right module structure $\psi_R$ which can be lifted to view $N_{D_2}^{D_{2m}}\mM$ as a right $N^{D_{2m}}_e \iota_e^* \mM$-module. We define this module structure via the induced action map
 \begin{equation} \label{eqn: hr bar left module}
 N_{D_2}^{D_{2m}}\mM\ \square\ N^{D_{2m}}_e \iota_e^* \mM \xrightarrow{\cong}
 N_{D_2}^{D_{2m}}(\mM\ \square\ N^{D_2}_e \iota_e^* \mM) \xrightarrow{N_{D_2}^{D_{2m}}(\psi_R)} 
 N_{D_2}^{D_{2m}}\mM, 
 \end{equation}
 where the isomorphism on the left arises from the fact that the norm is symmetric monoidal.

Let $\zeta=e^{2\pi i/2m}$. Note that $\zeta D_2 \zeta^{-1}$ is also an order 2 subgroup of $D_{2m}$. Since it is distinct from $D_2$ only up to a change of generator, we get an equivalence of categories, 
\[
c_\zeta: \mathit{Mack}^{D_2} \rightarrow \mathit{Mack}^{\zeta D_2 \zeta^{-1}}. 
\]
This equivalence may similarly be defined in spectra and we will use the same notation to denote it.

In the bar construction defining Real Hochschild homology, the right coefficients are given by the norm $\zetanorm \mM$. The given left module structure $\psi_L: N^{D_2}_e \iota_e^* \mM\ \square\ \mM \rightarrow \mM$ defines a left $N_e^{D_2} \iota_e^* \mM$-module structure on $\zetanorm \mM$ by composing the isomorphism,  
\begin{equation} \label{eqn: hr bar right module 1}
N^{D_2}_e \iota_e^* \mM\ \square\ \zetanorm \mM \cong N^{D_{2m}}_{\zeta D_2 \zeta^{-1}} (N^{\zeta D_2 \zeta^{-1}}_e \iota_e^* \mM\ \square\ c_{\zeta} \mM)
\end{equation}
with the induced left action,
\begin{equation}\label{eqn: hr bar right module 2}
N^{D_{2m}}_{\zeta D_2 \zeta^{-1}} (N^{\zeta D_2 \zeta^{-1}}_e \iota_e^* \mM\ \square\ c_{\zeta} \mM) \xrightarrow {N^{D_{2m}}_{\zeta D_2 \zeta^{-1}}(c_{\zeta} (\psi_L))} \zetanorm \mM.
\end{equation} 

With this description of the coefficients in the bar construction, we now recall the definition of Real Hochschild homology due to \cite{realhh}. 

\begin{defn} \label{defn: real hochschild homology}
Let $\mM$ be a discrete $E_{\sigma}$-ring. Let $\mhr_{\bullet}^{D_{2m}}(\mM)$ denote the two-sided bar construction 
\[
B_{\bullet}(N^{D_{2m}}_{D_2} \mM, N^{D_{2m}}_e \iota_e^* \mM, \zetanorm \mM)
\]
and define the \emph{Real $D_{2m}$-Hochschild homology of $\mM$} to be the $\Z$-graded $D_{2m}$-Mackey functor
\[
\mhr_*^{D_{2m}}(\mM) := H_*(\mhr_{\bullet}^{D_{2m}}(\mM)),
\]
where $H_*$ denotes taking the homology of the dg Mackey functor associated to $\mhr_{\bullet}^{D_{2m}}(\mM)$. 
\end{defn}

\begin{rem}
Theorem 6.20 of \cite{realhh} gives a Real linearization map 
    \[
    \mpi_n^{D_{2m}} \thr(A) \rightarrow \mhr_n^{D_{2m}}(\mpi_0^{D_2} A), 
    \]
which is one piece of evidence justifying that $\mhr$ is the proper algebraic analogue of $\thr$. The construction of the Real \bok\ spectral sequence in the following section relates $\mhr$ to the equivariant homology of $\thr$, further demonstrating that $\mhr$ is a suitable algebraic analogue of $\thr$. 
\end{rem}


\subsection{Real Hochschild homology for graded inputs}\label{section: graded inputs for hh}

The spectral sequence construction in the following section requires us to make sense of Real Hochschild homology for a graded discrete $E_{\sigma}$-ring. In the classical spectral sequence we only consider $\Z$-graded inputs; in the construction of equivariant \bok\ spectral sequences, the necessary flatness conditions that arise when one does not restrict to field coefficients are more likely to hold in an equivariant grading. For example, the equivariant \bok\ spectral sequence for $\thh_{C_n}$ constructed in \cite{twistedthh} takes inputs in $RO(C_n)$-graded Green functors.   

As we saw in Definition \ref{defn: real hochschild homology}, however, $\mhr$ is defined with a two-sided bar construction that involves taking equivariant norms and restrictions of the input. This change of groups complicates the question of gradings. Let $E$ be a commutative $D_{2m}$-ring spectrum and $A$ a ring spectrum with anti-involution. Our spectral sequence construction requires us to take the $E$-homology of the norm $N_{D_2}^{D_{2m}}A$ and identify it with the norm of the $\iota_{D_2}^*E$-homology of $A$. We must therefore extend the grading to a larger setting, $\underline{RO}(G)$.

Let $G$ be a finite group. An element in $\underline{RO}(G)$ is a pair $(H, \alpha)$ where $H \subseteq G$ and $\alpha$ is a virtual $H$-representation. Note the contrast with an element of $RO(G)$ which is only a virtual $G$-representation. Since $\underline{RO}(G)$ contains representations of all groups involved in the norm, it is the appropriate grading in this equivariant setting with norms. 

This expanded grading has been considered in \cite{hill2017slice}, \cite{rograded}, and \cite{hillfreeness}. For a detailed construction and definition of $\underline{RO}(G)$-graded Mackey functors, we direct the reader to work of Angeltveit-Bohmann in \cite{rograded}. We recall here the particular instance of the $\underline{RO}(G)$-graded homotopy Mackey functor.

We have the following definition of the $\underline{RO}(G)$-graded homotopy of a $G$-spectrum due to Hill-Hopkins-Ravenel \cite{hill2017slice}.

\begin{defn}\label{defn: RO graded homotopy}
    Let $X$ be a $G$-spectrum. For each pair $(H, \alpha)$ consisting of a subgroup $H \subseteq G$ and a virtual orthogonal representation $\alpha$ of $H$, define the \emph{$\underline{RO}(G)$-graded homotopy Mackey functor} to be the $G$-Mackey functor $\mpi_{H, \alpha}(X)$ where
    \[
    \mpi_{H, \alpha} (X)(T) := [(G_+ \smsh_H S^{\alpha}) \smsh T_+, X]_G \cong [S^{\alpha} \smsh \iota_H^* T_+, \iota_H^*X]_H=\mpi_{\alpha}^H(\iota_H^*X)(\iota_H^*T)
    \]
    for each finite $G$-set $T$. Taken together, this forms an $\underline{RO}(G)$-graded Mackey functor which is denoted by $\mpi_{\ostar}(X)$. 
\end{defn}
 
\begin{nota}
    We will use $\ostar$ to denote an $\underline{RO}(G)$-grading. We continue to denote an $RO(G)$-grading by $\star$.  
\end{nota}

\begin{defn} \label{defn: RO graded homology}
    For a $G$-spectrum $X$ and a commutative $G$-ring spectrum $E$, the \emph{$\underline{RO}(G)$-graded $E$-homology of $X$} is defined to be
    \[
    \underline{E}_{\ostar}(X) := \mpi_{\ostar}(X \smsh E).
    \]
\end{defn}

Recall that the category of $G$-Mackey functors is equipped with a symmetric monoidal product called the box product. Similarly, this category of $\underline{RO}(G)$-graded Mackey functors has a symmetric monoidal product (see, for instance \cite{hillfreeness}). For $\underline{RO}(G)$-graded Mackey functors $\mM_{\ostar}$ and $\mN_{\ostar}$ we denote this product, which we will refer to as the \emph{$\underline{RO}(G)$-graded box product}, by $\mM_{\ostar} \square\ \mN_{\ostar}$. The notion of a product allows to consider monoids in this category.  

\begin{defn} \label{defn: RO graded green functor}
    An \emph{$\underline{RO}(G)$-graded Green functor} is an associative monoid in the category of $\underline{RO}(G)$-graded Mackey functors with respect to the graded box product.  
\end{defn}

This definition gives rise to a notion of $\mR_{\ostar}$-modules for an $\underline{RO}(G)$-graded $G$-Green functor.

\begin{defn}
    Let $\mR_{\ostar}$ be an $\underline{RO}(G)$-graded $G$-Green functor. A \emph{left $\mR_{\ostar}$-module} is an $\underline{RO}(G)$-graded Mackey functor $\mM_{\ostar}$ with action map 
    \[ \mR_{\ostar} \square \mM_{\ostar} \rightarrow \mM_{\ostar} \]
    which satisfies the usual module relations but over the graded box product. A right module is defined analogously. 
 \end{defn}
 
 Let $\mN_{\ostar}$ be a right $\mR_{\ostar}$-module with action $\psi$ and $\mM_{\ostar}$ a left module with action $\phi$. The relative graded box product  $\mN_{\ostar} \boxover{\mR_{\ostar}} \mM_{\ostar}$ is the coequalizer

\begin{center}
\begin{tikzcd}
\mN_{\ostar} \square \mR_{\ostar} \square \mM_{\ostar} \arrow[r, shift left=1, "\psi \square id"] \arrow[r, shift right=1, "id \square \phi"'] & \mN_{\ostar} \square \mM_{\ostar} \arrow[r, rightarrow] & \mN_{\ostar} \square_{\mR_{\ostar}} \mM_{\ostar}. 
\end{tikzcd}
\end{center}
  
 \begin{defn} \label{defn: RO graded flat module}  
    We say that an $\underline{RO}(G)$-graded $\mR_{\ostar}$-module $\mM_{\ostar}$ is \emph{flat} if the functor $(-) \boxover{\mR_{\ostar}} \mM_{\ostar}$ is exact.
\end{defn}

Our motivation for the use of these $\underline{RO}(G)$-graded objects is to define an appropriately graded input to Real Hochschild homology, the construction of which includes equivariant norms. We now extend the work \cite{realhh} and define a notion of $\underline{RO}(D_2)$-graded discrete $E_{\sigma}$-rings. 

\begin{defn} \label{defn: ro graded discrete e sigma ring}
    An \emph{$\underline{RO}(D_2)$-graded discrete $E_{\sigma}$-ring} $\mM_{\ostar}$ is the following data:
    \begin{enumerate}
        \item A $D_2$-Mackey functor $\mM_{(H, \alpha)}$ for each subgroup $H \subseteq D_2$ and virtual $H$-representation $\alpha$ such that $\mM_{\ostar}(D_2/e)$ forms an $\underline{RO}(D_2)$-graded ring with anti-involution. That is, we have an associative product, 
        \[
        \mM_{\ostar}(D_2/e) \otimes \mM_{\ostar}(D_2/e) \rightarrow \mM_{\ostar}(D_2/e)
        \]
        where the domain has the action of swapping the two copies of $\mM_{\ostar}(D_2/e)$. 
        \item An $N^{D_2}_e \iota_e^* \mM_{\ostar}$-bimodule structure on $\mM_{\ostar}$. We further require that the action restricts to the standard action of $\mM_{\ostar}(D_2/e) \otimes \mM_{\ostar}(D_2/e)^{op}$ on $\mM_{\ostar}(D_2/e)$. 
        \item A designated unit $1 \in \mM_{\ostar}(D_2/D_2)$ which restricts to $1 \in \mM_{\ostar}(D_2/e)$.  
    \end{enumerate}
\end{defn}

We claim that the $\underline{RO}(D_2)$-graded homotopy Mackey functor of a ring spectrum with anti-involution forms an $\underline{RO}(D_2)$-graded discrete $E_{\sigma}$-ring. A useful perspective in proving this statement is that ring spectra with anti-involution are algebras in $D_2$-spectra over an $E_{\sigma}$-operad. It is this interpretation which motivates the name and definition of a discrete $E_{\sigma}$-ring in equivariant algebra. Formally, this gives the following definition:

\begin{defn}[\cite{realhh}, Corollary 3.10] \label{defn: e sigma ring}
   An \emph{$E_{\sigma}$-ring} $A$ is a $D_2$-spectrum such that 
   \begin{enumerate}
       \item the spectrum $\iota_e^*A$ is an $E_1$-ring with anti-involution, denoted $\tau: \iota_e^* A^{op} \rightarrow \iota_e^* A$ and given by the action of the generator of the Weyl group.
       \item the spectrum $A$ is an $E_0-N^{D_2}_e\iota_e^*A$-algebra and applying the restriction functor $\iota_e^*$ to the $N^{D_2}_e\iota_e^*A$-module structure map gives $\iota_e^*A$ the standard $\iota_e^*A \smsh \iota_e^* A^{op}$-module structure.  
   \end{enumerate}
\end{defn}

\begin{prop}\label{thm: homotopy of rswai is discrete e sigma}
    The $\underline{RO}(D_2)$-graded homotopy of a ring spectrum with anti-involution (equivalently, an $E_{\sigma}$-ring) $A$ forms an $\underline{RO}(D_2)$-graded discrete $E_{\sigma}$-ring. 
\end{prop} 
 
\begin{proof}
    By Definition \ref{defn: RO graded homotopy} we see that each $\mpi_{(H, \alpha)}(A)$ is a $D_2$-Mackey functor. Further, we see by this definition that the restriction to the orbit $(D_2/e)$ recovers the non-equivariant homotopy of $\iota_e^* A$. This is clear in the case of $H=e$. When $H=D_2$ we have that 
    \[
    \mpi_{(D_2, \alpha)}(A)(D_2/e) = [S^{\alpha}, A]_e = \pi_{|\alpha|}(\iota_e^*A). 
    \]
    Thus, restricting to the orbit $(D_2/e)$ in both cases recovers the non-equivariant homotopy of the underlying spectrum. From Definition \ref{defn: e sigma ring}, we know that $\iota_e^* A$ is an $E^1$-ring with anti-involution. Thus we have a map $\iota_e^* A \smsh \iota_e^* A \rightarrow A$ with a swap action on the domain. This induces the desired map on homotopy. 
    
    To define the module structure specified in condition (2) of Definition \ref{defn: ro graded discrete e sigma ring}, we recall that since $A$ is an $E_{\sigma}$-ring, it has an $N^{D_2}_e \iota_e^* A$-bimodule structure. We denote these module action maps by $\psi_R'$ and $\psi_L'$. These maps induce module action maps on $\underline{RO}(D_2)$-graded homotopy,
     \[
     \psi_R': {\mpi}_{\ostar}(A \smsh N^{D_2}_e \iota_e^* A) \rightarrow {\mpi}_{\ostar}(A),
     \]
     and similarly for $\psi_L'$. Since the $\underline{RO}(D_2)$-graded homotopy functor is $D_2$-lax monoidal there is a map
     \[
        {\mpi}_{\ostar}(A) \square\ N^{D_2}_e \iota_e^* {\mpi}_{\ostar}(A) \rightarrow {\mpi}_{\ostar}(A \smsh N^{D_2}_e \iota_e^* A). 
     \]
     Postcomposing this map with $\psi_R'$ yields the desired right-module structure,
     \[
     \psi_R: {\mpi}_{\ostar}(A) \square\ N^{D_2}_e \iota_e^*{\mpi}_{\ostar}(A) \rightarrow {\mpi}_{\ostar}(A). 
     \]
     The left module action $\psi_L$ is defined analogously. Thus ${\mpi}_{\ostar}(A)$ is an $N^{D_2}_e \iota_e^* {\mpi}_{\ostar}(A)$-bimodule. Above, we argued that restriction to the orbit $(D_2/e)$ recovers the non-equivariant homotopy of $\iota_e^* A$. By Definition \ref{defn: e sigma ring}, the restriction of the module structure map recovers the standard action. Thus the induced map on homotopy is also the standard action of $\mpi_{\ostar}(A)(D_2/e) \otimes \mpi_{\ostar}(A)(D_2/e)^{op}$. 
\end{proof}

The $\underline{RO}(D_2)$-graded equivariant homology of a ring spectrum with anti-involution will be the input in the Real \bok\ spectral sequence. We now define $\mhr$ for $\underline{RO}(D_2)$-graded discrete $E_{\sigma}$-rings. We use the $N^{D_2}_e \mpi_{\ostar}(A)$-module structure of $\mpi_{\ostar}(A)$ to give the necessary module structures for the coefficients in the two-sided bar construction. 

\begin{prop}
    Let $\mM_{\ostar}$ be an $\underline{RO}(D_2)$-graded discrete $E_{\sigma}$-ring. Then $N^{D_{2m}}_{D_2}\mM_{\ostar}$ is a right $N^{D_{2m}}_e \iota_e^* \mM_{\ostar}$-module and $\zetanorm \mM_{\ostar}$ is a left $N^{D_{2m}}_e \iota_e^* \mM_{\ostar}$-module.
\end{prop}

\begin{proof}
    Since $\mM_{\ostar}$ is an $\underline{RO}(D_2)$-graded discrete $E_{\sigma}$-ring there is an $N^{D_2}_e \iota_e^* \mM_{\ostar}$-bimodule structure on $\mM_{\ostar}$, 
    \begin{align*} 
    \psi_R'&: N^{D_2}_e \iota_e^* \mM_{\ostar} \square\ \mM_{\ostar} \rightarrow \mM_{\ostar}\\
    \psi_L' &: \mM_{\ostar} \square\ N^{D_2}_e \iota_e^* \mM_{\ostar} \rightarrow \mM_{\ostar}.  
   \end{align*}
    We then define the desired module structures over the $\underline{RO}(D_2)$-graded box product to be the same composites as specified in the ungraded case (see Equations \ref{eqn: hr bar left module}, \ref{eqn: hr bar right module 1}, and \ref{eqn: hr bar right module 2}).    
\end{proof}

\begin{defn}\label{defn: graded hr}
For an $\underline{RO}(D_2)$-graded discrete $E_{\sigma}$-ring $\mM_{\ostar}$, we define the two-sided bar construction
\[
\mhr^{D_{2m}}_{\bullet}(\mM_{\ostar}):= B_{\bullet}(N^{D_{2m}}_{D_2} \mM_{\ostar}, N^{D_{2m}}_e \iota_e^* \mM_{\ostar}, \zetanorm \mM_{\ostar})
\]
with the same face and degeneracy maps as in the ungraded cases, taken here over the $\underline{RO}(D_2)$-graded box product. The \emph{Real Hochschild homology of the $\underline{RO}(D_2)$-graded discrete $E_{\sigma}$-ring} $\mM_{\ostar}$ is the homology of this two-sided bar construction. 
\end{defn}


\section{Spectral sequence construction}\label{section: construction}

We are now ready to construct a \bok\ spectral sequence for Real topological Hochschild homology. The presence of equivariant norms in these constructions are notable complications arising in the arguments which otherwise follow the original ones of \bok. These norms require us to place additional hypotheses on the inputs to ensure that we may recognize the $E^2$-page of the spectral sequence as Real Hochschild homology. The existence of norms on the $E^2$-page also necessitates the use of $\underline{RO}$-graded homology theories. We begin this section by recalling Hill's notion of free homology in an equivariant setting, particularly as it relates to this question of how to treat equivariant norms in the spectral sequence construction. 


\subsection{Free homology}

\begin{defn}\label{defn: free homology}
    Let $G$ be a finite group, $A$ be a $G$-spectrum, and $E$ be a commutative $G$-ring spectrum. We say $A$ has \textit{free $E$-homology} if $E \smsh A$ splits as a wedge of $E$-modules of the form
    \[
    E \smsh (G_+ \smsh_H S^{\alpha})
    \]
    where $\alpha$ is a virtual representation of $H$, a subgroup of $G$.
\end{defn}
Hill shows that this class of spectra with $E$-free homology is particularly nice, enjoying the property of closure under norms. Of particular relevance to our work is a description of how the homology functor interacts with the equivariant norm under freeness hypotheses. 

\begin{lemma}[\cite{hillfreeness}, Corollary 3.30] \label{thm: norm commutes w homology}
  Let $A$ be an $H$-spectrum for $H \subseteq G$ and let $E$ be a commutative $G$-ring spectrum. If $A$ has free $(\iota_H^*E)$-homology, then there is a natural isomorphism
  \[
  \underline{E}_{\ostar}(N^G_H A) \cong N^G_H(\underline{(\iota_H^* E)}_{\ostar} (A)). 
  \]
\end{lemma}

This result, which allows us to permute the norm functor and the homology functor, will be important in our construction of the Real \bok\ spectral sequence. In particular, we make use of the following consequence of this lemma. 

\begin{cor} \label{thm: trivial homology is always free}
    Let $A$ be a ring spectrum with anti-involution and $E$ be a commutative $D_{2m}$-ring spectrum such that $\iota_e^*A$ has free $\iota_e^*E$-homology. Then there is an isomorphism
    \[
    \underline{E}_{\ostar}(N_e^{D_{2m}}\iota_e^* A) \cong N_e^{D_{2m}}\iota_e^* (\underline{(\iota_{D_2}^*E)}_{\ostar} (A))  
    \]
\end{cor}

\begin{proof}
    The assumption that $\iota_e^*A$ has free $\iota_e^*E$-homology allows us to apply the result in Lemma \ref{thm: norm commutes w homology}. We have
    \[
    \underline{E}_{\ostar}(N_e^{D_{2m}}\iota_e^* A) \cong N_e^{D_{2m}}(\underline{\iota_e^* E}_{\ostar} (\iota_e^*A)).
    \]
    The restriction functor commutes with homology so on the right hand side we can write
    \[
    N_e^{D_{2m}}(\underline{\iota_e^* E}_{\ostar} (\iota_e^*A)) \cong N_e^{D_{2m}}\iota_e^* (\underline{(\iota_{D_2}^*E)}_{\ostar} (A)).
    \]
    which gives us the desired isomorphism. 
\end{proof}

\begin{rem} \label{rem: always free HFp}
This freeness condition is always satisfied when $E=H\mF_p$ since $\iota_e^* H\mF_p \smsh A$ is considered by Ravenel in \cite{ravenel2023complex} Proposition 2.1.2.g and shown to be a wedge of suspensions of $\iota_e^* H\mF_p$.
\end{rem}


\subsection{Construction of the Real \bok\ spectral sequence}
We first recall the classical \bok\ spectral sequence, which was developed by \bok\ in \cite{bokstedtthhZandZp} and used to compute $\thh$ of $\mathbb{F}_p$ and $\mathbb{Z}$. For a more modern treatment of the \bok\ spectral sequence construction using the language of spectra and smash products, see \cite{ekmm}. 

\begin{thm}[\cite{ekmm}, Theorem IX.2.9]\label{thm: classical bok ss}
Let $E$ be a commutative ring spectrum, $A$ a ring spectrum, and $M$ a cellular $A$-bimodule. If $E_*(A)$ is $E_*$-flat, then there is a spectral sequence of the form
\[
E^2_{p,q} = \hh^{E_*}_{p,q} (E_*(A); E_*(M)) \Rightarrow E_{p+q}(\thh(A;M)).
\]
\end{thm}

The classical \bok\ spectral sequence arises from a more general spectral sequence which computes the $E$-homology of a proper simplicial spectrum. A simplicial spectrum $X_{\bullet}$ is called \emph{proper} if the inclusion of the degenerate subspectrum $sX_q \rightarrow X_q$ is a cofibration for each $q \geq 0$ (see \cite{ekmm}, X.2.1). 

\begin{thm}[\cite{ekmm}, Theorem X.2.9]\label{thm: ekmm simp ss}
Let $X_{\bullet}$ be a proper simplicial spectrum and let $E$ be any spectrum. Then there is a natural homological spectral sequence 
\[
E^2_{p,q} = H_p(E_q(X_{\bullet})) \Rightarrow E_*(|X_{\bullet}|).
\]
\end{thm}

We now proceed to construct a Real equivariant version of the \bok\ spectral sequence. To begin, we employ a result of Angelini-Knoll, Gerhardt, and Hill to describe $\thr$ via the multiplicative double coset formula. 

\begin{lemma}[\cite{realhh}, Theorem 5.9] \label{thm: mult double coset thr}
    Let $A$ be a flat ring spectrum with anti-involution. There is a stable equivalence of $D_{2m}$-spectra
    \[
   \iota_{D_{2m}}^*N^{O(2)}_{D_2}A \simeq N^{D_{2m}}_{D_2}A \smsh^{\mathbb{L}}_{N^{D_{2m}}_e \iota_e^*A} \zetanorm A.  
    \]
    In particular, this gives an equivalence
    \[
    \iota_{D_{2m}}^*\thr(A) \simeq |B_{\bullet}(N^{D_{2m}}_{D_2}A, N^{D_{2m}}_e \iota_e^* A, \zetanorm A)|.
    \]
\end{lemma}

\begin{rem}
 In Proposition \ref{thm: homotopy of rswai is discrete e sigma} we showed that the $\underline{RO}(D_2)$-graded homotopy of a ring spectrum with anti-involution $A$ forms an $\underline{RO}(D_2)$-graded discrete $E_{\sigma}$-ring. Note that for $E$ a commutative $D_{2m}$-ring spectrum, the $\underline{RO}(D_{2m})$-graded $E$-homology of $A$ is also an $\underline{RO}(D_2)$-graded discrete $E_{\sigma}$-ring. This is because $E$-homology is defined by taking the $\underline{RO}(D_2)$-graded homotopy of the spectrum $A \smsh E$ which is also a ring spectrum with anti-involution.  
\end{rem}

\begin{nota}
    For a $G$-spectrum $E$, we let $\underline{E}_{\ostar}$ denote the $\underline{RO}(G)$-graded equivariant homotopy Mackey functor of $E$, $\mpi^G_{\ostar}(E)$. 
\end{nota}

One may consider the symmetric monoidal category of $\underline{E}_{\ostar}$-algebras. The product in this category is a box product over $\underline{E}_{\ostar}$, denoted $\square_{\underline{E}_{\ostar}}$, which is defined analogously to the relative tensor product. When an $\underline{RO}(D_2)$-graded discrete $E_{\sigma}$-ring $\underline{A}_{\ostar}$ is an $\underline{E}_{\ostar}$-algebra, we can define Real Hochschild homology $\underline{\hr}_*^{\underline{E}_{\ostar}, D_{2m}}(\underline{A}_{\ostar})$ as in Definition \ref{defn: graded hr} using a graded box product over $\underline{E}_{\ostar}$.

The main result of this paper is the following theorem. 

\begin{thm} \label{thm: real bok ss}
    Let $A$ be a ring spectrum with anti-involution and let $E$ be a commutative $D_{2m}$-ring spectrum. If $\underline{E}_{\ostar}(N^{D_{2m}}_{D_2} A)$ and $\underline{E}_{\ostar}(N^{D_{2m}}_e \iota^*_e A)$ are both flat as modules over $\underline{E}_{\ostar}$ and if $A$ has free $(\iota_{D_2}^*E)$- and $\iota_{e}^*E$-homology then there is a Real \bok\ spectral sequence of the form,
    \[
    E^2_{*, \ostar} = \mhr^{\underline{E}_{\ostar}, D_{2m}}_*(\underline{(\iota_{D_2}^*E)}_{\ostar}(A)) \Rightarrow \underline{E}_{\ostar}(\iota_{D_{2m}}^* \thr(A))
    \]
    where $\underline{E}_{\ostar}$ denotes the $\underline{RO}(D_{2m})$-graded homotopy Mackey functor $\mpi_{\ostar}(E)$.
\end{thm}

\begin{proof} Note that the proof of Theorem \ref{thm: ekmm simp ss} goes through equivariantly. Therefore, if $E$ is a commutative $G$-spectrum and $X_{\bullet}$ is a proper simplicial $G$-spectrum, there is a spectral sequence 
\[
E^2_{*,\ostar}=H_*(\underline{E}_{\ostar}(X_{\bullet})) \Rightarrow E_{*+\ostar}(|X_{\bullet}|). 
\]  

As in the non-equivariant case, the simplicial filtration
\[
...\ F_{p-1} \subset F_{p} \subset F_{p+1} \subset ... \subset X_{\bullet}
\]
gives rise to such a spectral sequence where the $E^1$-page is given by
\[
E^1_{p, \ostar} = \underline{E}_{\ostar}(F_{p}/F_{p-1}). 
\]

Since the multiplicative double coset formula in Lemma \ref{thm: mult double coset thr} is a simplicial spectrum, we can apply this result to construct the Real \bok\ spectral sequence. 

    We begin by taking $E$-homology of the double bar construction from the multiplicative double coset formula for $\iota_{D_{2m}}^* \thr(A)$. At the $p$th level we see that $E^1_{p, \ostar}$ takes the form:
\begin{align*}
    & \underline{E}_{\ostar}(N^{D_{2m}}_{D_2}A \smsh \overbrace{N^{D_{2m}}_e \iota_e^* A \smsh ... \smsh N^{D_{2m}}_e \iota_e^* A}^{p} \smsh \zetanorm A)\\
   \cong\ & \mpi_{\ostar}( N^{D_{2m}}_{D_2}A \smsh {E} \smsh_{{E}} N^{D_{2m}}_e \iota_e^* A \smsh E \smsh_{E} ... \smsh_{E} N^{D_{2m}}_e \iota_e^* A \smsh E \smsh_{E} \zetanorm A \smsh E) \\
  \cong\ & \mpi_{\ostar}(N^{D_{2m}}_{D_2}A \smsh E) \square \mpi_{\ostar}(N^{D_{2m}}_e \iota_e^* A \smsh E) \square ...
    \square \mpi_{\ostar}(N^{D_{2m}}_e \iota_e^* A \smsh E) \square \mpi_{\ostar}(\zetanorm A \smsh E) \\
   =\ &\underline{E}_{\ostar}(N^{D_{2m}}_{D_2}A) \square \underline{E}_{\ostar}(N^{D_{2m}}_e \iota_e^* A) \square ... \square \underline{E}_{\ostar}(N^{D_{2m}}_e \iota_e^* A) \square \underline{E}_{\ostar}(\zetanorm A).   
\end{align*}

All box products above are taken over $\underline{E}_{\ostar}$; we omit this from the notation for the ease of reading. 

By Corollary 3.20 of \cite{hillfreeness} and the flatness of $\underline{E}_{\ostar}(N_{D_2}^{D_{2m}}A)$ and $\underline{E}_{\ostar} (N_e^{D_{2m}} \iota_e^* A),$ the K{\"u}nneth isomorphism in \cite{lewismandell} gives rise to the second isomorphism above. 

Recall from Corollary \ref{thm: trivial homology is always free} that when $A$ has free $\iota_e^* E$-homology as we assumed, there is an isomorphism
\begin{equation} \label{eqn: norm 1}
   \underline{E}_{\ostar}(N_e^{D_{2m}}\iota_e^* A) \cong N_e^{D_{2m}}\iota_e^* (\underline{(\iota_{D_2}^*E)}_{\ostar} (A)).
\end{equation}
Furthermore, the hypothesis that $A$ has free $(\iota_{D_2}^*E)$-homology in conjunction with Lemma \ref{thm: norm commutes w homology}, yields an isomorphism 
\begin{equation} \label{eqn: norm 2}
      \underline{E}_{\ostar}(N^{D_{2m}}_{D_2}A) \cong N^{D_{2m}}_{D_2} (\underline{(\iota_{D_2}^*E)}_{\ostar} (A)).
\end{equation}
Finally, the equivalence of categories defined by $c_{\zeta}$ implies that $c_{\zeta} A$ has free $(\iota_{\zeta D_2 \zeta^{-1}}^*E)$-homology. Thus we have an isomorphism,
\begin{equation} \label{eqn: norm 3}
   \underline{E}_{\ostar}(\zetanorm A) \cong \zetanorm (\underline{(\iota_{D_2}^*E)}_{\ostar} (A)).
\end{equation}

The isomorphisms in \ref{eqn: norm 1}, \ref{eqn: norm 2}, and \ref{eqn: norm 3} allow us to conclude that $E^1_{p, \ostar}$ is isomorphic to 
\[
N^{D_{2m}}_{D_2} (\underline{(\iota_{D_2}^*E)}_{\ostar} (A)) \square  
{N_e^{D_{2m}}\iota_e^* (\underline{(\iota_{D_2}^*E)}_{\ostar} (A))}^{\square  p }
\square 
\zetanorm (\underline{(\iota_{D_2}^*E)}_{\ostar} (A)),
\]
where again, the box product occurs over $\underline{E}_{\ostar}$.
This is precisely the $p$th simplicial level of the two-sided bar construction which defines the Real $D_{2m}$-Hochschild homology,
\[
\mhr^{D_{2m}}_{\bullet}(\underline{(\iota_{D_2}^*E)}_{\ostar} (A)).
\] 
At each simplicial level, we have identified the $E_1$-page of the spectral sequence with the complex that computed Real Hochschild homology. A diagram chase shows that the $d_1$ differential of the spectral sequence agrees with the differential in the complex for $\mhr$ and hence on the $E_2$-page we have Real Hochschild homology.   
\end{proof}

A particular case of interest occurs when $E=H\mF_2$ and $m=1$; the Real \bok\ spectral sequence allows us to compute the $D_2$-equivariant homology of $\thr(A)$ as a $D_2$-spectrum. Note that in this case we do not need an additional flatness assumption about $\underline{E}_{\ostar}(N^{D_2}_{D_2} A)$ since this norm is trivial. Similarly, we do not require freeness assumptions since we are in the case described by Remark \ref{rem: always free HFp}. Thus we have the following corollary:

\begin{cor}
    Let $A$ be a ring spectrum with anti-involution and such that $\underline{H}^{D_2}_{\ostar}(A; \mF_2)$ and $\underline{H}^{D_2}_{\ostar}(N^{D_2}_e \iota_e^* A; \mF_2)$ are flat as modules over ${H\mF_2}_{\ostar}$. Then there is a Real \bok\ spectral sequence
    \[
    E^2_{*, \ostar}=\mhr^{D_2}_*(\underline{H}_{\ostar}^{D_2} (A; \mF_2)) \Rightarrow \underline{H}_{\ostar}^{D_2}(\iota_{D_2}^* \thr(A); \mF_2).
    \]
\end{cor}

\begin{rem} 
The conditions of flatness and freeness in Theorem \ref{thm: real bok ss} are not so restrictive that we cannot consider interesting inputs. One example is the Real bordism spectrum, $\bord$, a $D_2$-equivariant ring spectrum first studied in \cite{landweber1968conjugations} and \cite{fujii1976cobordism}, which played an important role in work of Hill-Hopkins-Ravenel on the Kervaire invariant one problem \cite{hhrkervaire}. Since $\underline{H}^{D_2}_{\ostar}(\bord; \mF_2)$ is polynomial over $H\mF_2$, it satisfies the flatness condition required to use the Real \bok\ spectral sequence.
\end{rem}



\section{Extensions to Twisted \texorpdfstring{$\thh$}{THH}} \label{section: twisted thh}
The techniques used in the previous section to construct the Real \bok\ spectral sequence can be extended to another equivariant generalization of topological Hochschild homology which takes inputs in ring spectra with a $C_n$-action. In \cite{twistedthh} the authors construct a \bok\ spectra sequence for this $C_n$-twisted $\thh$. In this section, we generalize this spectral sequence to one which computes the $G$-equivariant homology of $H$-twisted $\thh$ for $H \subseteq G$ subgroups of $S^1$. 

\subsection{Twisted topological Hochschild homology}

We begin by recalling the theory of $C_n$-twisted topological Hochschild homology, denoted by $\thh_{C_n}$, first defined in \cite{tcvianorm}.

\begin{defn}\label{defn: twisted cylic bar spectra}
    Let $R$ be an associative orthogonal $C_n$-ring spectrum. The \emph{$C_n$-twisted cyclic bar construction} on $R$, denoted by $B^{cy, C_n}_{\bullet}(R)$, is a simplicial spectrum which has $k$-simplices
    \[
    B^{cy, C_n}_{k}(R)= R^{\smsh (k+1)}. 
    \]
  For $g$ a generator of $C_n$, let $\alpha_k$ denote the composition that wraps the last factor of $R$ around to the front then acts by $g$: 
   \[
    \alpha_k: R \smsh \overbrace{R \smsh ... \smsh R}^{k} \rightarrow R \smsh R \smsh \overbrace{R \smsh ... \smsh R}^{k-1} \rightarrow {}^g R \smsh R \smsh \overbrace{R \smsh ... \smsh R}^{k-1}.
   \]
   We define the face maps of $B^{cy, C_n}_{k}(R)$ by
   \begin{align*}
    d_i= \begin{cases} 
     id^{\smsh i} \smsh \mu \smsh id^{\smsh (k-i-1)} & 0\leq i<k \\
      (\mu \smsh id^{\smsh (k-1)}) \circ \alpha_k & i=k. 
   \end{cases}
\end{align*}
The degeneracy maps are defined as usual.
\end{defn}

\begin{defn}\label{defn: twisted thh}
  Let $R$ be an associative orthogonal $C_n$-ring spectrum and define the \emph{$C_n$-twisted topological Hochschild homology} of $R$ to be     
\[
\thh_{C_n}(R) = N^{S^1}_{C_n}(R)= |B^{cy, C_n}_{\bullet}|.
\]
\end{defn}

There is an algebraic theory associated to $\thh_{C_n}$, which takes inputs in $C_n$-Green functors. This theory of Hochschild homology for Green functors was first defined by Blumberg-Gerhardt-Hill-Lawson \cite{hhforgfs} and we now recall some important definitions. 

 Let $G  \subset S^1$ be a finite subgroup and let $g \in G$. For a $G$-Green functor $\mR$ and a left module $\psi: \mR \Box \mM \rightarrow \mM$, we can define a $g$-twisted module structure on $\mM$, denoted ${}^g\mM$, where the action map ${}^g\psi$ is the composition
    \begin{center}
    \begin{tikzcd}
        \mR \Box \mM \arrow[d, "g \Box 1"'] \arrow[dr, "{}^g\psi"] &  \\
        \mR \Box \mM \arrow[r, "\psi"] & \mM.
    \end{tikzcd}
    \end{center}

\begin{defn} \label{defn: hh for green functors}
   Let $G \subset S^1$ be a finite subgroup and let $g=e^{2 \pi i/|G|}$ be a generator of $G$. For a $G$-Green functor $\mR$ and an $\mR$-bimodule $\mM$, we denote the \emph{$G$-twisted cyclic nerve} by $\underline{B}_{\bullet}^{cy,G}(\mR;{}^g\mM)$. This is a simplicial Mackey functor which has $k$ simplices
   \[
   \underline{B}_k^{cy,G}(\mR;{}^g\mM)={}^g\mM \Box \mR^{\Box k}.   
   \]
The face map $d_0$ is given by the right module action. At simplicial level $k$, the face maps $d_i$ for $0< i <k$ are given by multiplying the $i$ and $(i+1)$st factors. The final face map $d_k$ rotates the $k$ factor to the $0$th position and then uses the $g$-twisted left module action. The degeneracy maps are the usual ones.  
\end{defn}

We now recall a relative version of $G$-twisted Hochschild homology where the input is an $H$-spectrum for $H \subseteq G$. 

\begin{defn} \label{defn: relative hh for green functors }
  Let $H \subseteq G$ be a finite subgroup of $S^1$ and let $g=e^{2 \pi i/|G|}$ be a generator of $G$ as above. For an associative $H$-Green functor $\mR$, we define the \emph{$G$-twisted Hochschild homology} of $\mR$ as the homology of the simplicial Mackey functor
  \[
  \mhh^G_H(\mR)_*=H_*(\underline{B}_{\bullet}^{cy,G}(N^G_H\mR; {}^gN^G_H\mR)).
  \] 
\end{defn}

\subsection{\texorpdfstring{$\underline{RO}(G)$}{RO(G)}-graded twisted \texorpdfstring{$\thh$}{THH}}

In the spectral sequence construction which follows, we again must consider a notion of Hochschild homology in the $\underline{RO}(G)$-graded setting. In the twisted cyclic bar construction which defines twisted $\thh$, the final face map picks up an additional sign in the graded case - see Section 4.1 of \cite{twistedthh} for a discussion of this point in the $\Z$- and $RO(G)$-graded settings. Before defining a notion of $\underline{RO}(G)$-graded Hochschild homology for $G$-Green functors, we address the question of signs in the final face map of the twisted cyclic bar construction in the $\underline{RO}(G)$-graded case. 

Let $G$ be a finite subgroup of $S^1$ and let $\beta$ and $\gamma$ be two finite dimensional real representations of $G$. The switch map on the representation spheres, 
\[S^{\beta} \smsh S^{\gamma} \rightarrow S^{\gamma} \smsh S^{\beta}\]
specifies an element in the Burnside ring $A(G)$, which we denote by $\sigma (\beta, \gamma)$. The rotating isomorphism of $\underline{RO}(G)$-graded Mackey functors is a map 
\[
\tau: \mM_{\ostar} \square\ \mN_{\ostar} \rightarrow \mN_{\ostar} \square\ \mM_{\ostar}.
\]
We restrict to working one subgroup at a time in the $\underline{RO}(G)$-graded box product, so at level $(H, \alpha)$ the rotating isomorphism is the same as in Definition 4.1.4 in \cite{twistedthh}, with the sign in the switch map coming from the Burnside ring $A(H)$. 

\begin{defn}
    We say an $\underline{RO}(G)$-graded Green functor $\mR_{\ostar}$ is \emph{commutative} if $\mu \tau = \mu$ where $\mu$ is the multiplication on $\mR_{\ostar}$ and $\tau$ is the rotating isomorphism $\mR_{\ostar} \square\ \mR_{\ostar} \rightarrow \mR_{\ostar}$. 
\end{defn}

With this description of how $\underline{RO}(G)$-graded Green functors commute past each other, we are ready to define an $\underline{RO}(G)$-graded notion of the twisted cyclic nerve. 

\begin{defn} \label{defn: ro graded relative hh for gfs}
    Let $G$ be a finite subgroup of $S^1$ and let $g=e^{2\pi i/|G|} \in G$ be a generator. For $\mR_{\ostar}$ an $\underline{RO}(G)$-graded $G$-Green functor and $\mM_{\ostar}$ an $\mR_{\ostar}$-module, we define the \emph{$G$-twisted cyclic nerve of $\mR_{\ostar}$ with coefficients in ${}^g\mM_{\ostar}$} to be the simplicial $\underline{RO}(G)$-graded Mackey functor which has $k$-simplices
    \[
    B^{cy, G}_{k}(\mR_{\ostar}, {}^g\mM_{\ostar})= {}^g\mM_{\ostar} \square\ \mR_{\ostar}^{\square k}. 
    \]
  The face map $d_0$ applies the right module action of $\mR_{\ostar}$ to ${}^g\mM_{\ostar}$. For $1 \leq i \leq k$, the face map $d_i$ multiplies the $i$th and $(i+1)$st copies of $\mR_{\ostar}$. The final face map $d_k$ incorporates the rotating isomorphism by rotating the last factor around to the front and then applying the left module action of $\mR_{\ostar}$ to ${}^g\mM_{\ostar}$. Explicitly, this is given by
\[
{}^g\mM_{\ostar} \square\ \mR_{\ostar}^{\square k} \xrightarrow{\tau_k} \mR_{\ostar} \square\ {}^g\mM_{\ostar} \square\ \mR_{\ostar}^{\square (k-1)} \xrightarrow{{}^g\psi \square id} {}^g\mM_{\ostar} \square\ \mR_{\ostar}^{\square (k-1)}.
\]
where $\tau_k$ denotes iterating the rotating isomorphism $k$ times in order to bring the last factor to the front and ${}^g \psi$ denotes the $g$-twisted module action on ${}^g\mM_{\ostar}$. The degeneracy maps in this simplicial object are the usual maps induced by inclusion via the unit. 
\end{defn}

In Section \ref{section: higher twisted bok ss}, we consider the case of relative Hochschild homology for $\underline{RO}(G)$-graded Green functors. With our definition of the graded twisted cyclic nerve, we can now define a relative equivariant Hochschild homology for these graded inputs. 

\begin{defn} \label{defn: hh for ro graded green functors}
    Let $H \subseteq G$ be finite subgroups of $S^1$ and let $\mR_{\ostar}$ be an $\underline{RO}(G)$-graded associative Green functor for $H$. The \emph{$G$-twisted Hochschild homology of $\mR_{\ostar}$} is 
    \[
    \mhh_H^G(\mR_{\ostar})_* := H_*(B^{cy, G}_{\bullet}(N^G_H \mR_{\ostar})).
    \]
\end{defn}

 \subsection{A \bok\ Spectral Sequence for Twisted \texorpdfstring{$\thh$}{THH}}\label{section: higher twisted bok ss}
In \cite{twistedthh}, the authors construct the following \bok\ spectral sequence for $G$-twisted topological Hochschild homology:

\begin{lemma}[\cite{twistedthh}, Theorem 4.2.7]\label{twistedbokss}
Let $G \subset S^1$ be a finite subgroup and $g=e^{2\pi i/|G|}$ a generator of $G$. Let $R$ be a $G$-ring spectrum and $E$ a commutative $G$-ring spectrum such that $g$ acts trivially on $E$. If $\underline{E}_{\star}(R)$ is flat over $\underline{E}_{\star}$, then there is a $G$-twisted \bok\ spectral sequence 
\[
E_{*,\star}^2=\mhh_*^{\underline{E}_{\star},G}(\underline{E}_{\star}(R)) \Rightarrow \underline{E}_{*+\star}(\iota^*_G\thh_G(R)).
\]
\end{lemma}

\begin{rem} \label{rem: trivial g action}

If $g$ acts trivially on $E$, Lemma 4.2.5 of \cite{twistedthh} shows that there is an isomorphism of left $\underline{E}_{\star}(R)$-modules, ${}^g\underline{E}_{\star}(R) \cong \underline{E}_{\star}({}^gR)$. The proof of this lemma gives an isomorphism at the level of spectra, $E \smsh {}^g R \cong {}^g(E \smsh R)$, hence the result also goes through when passing to $\underline{RO}(G)$-graded homotopy. 

 \end{rem}
  
This spectral sequence computes the equivariant homology of the $G$-restriction of $\thh_G$. However, since twisted $\thh$ is an $S^1$-spectrum, we could also consider the $G$-restriction of $\thh_H$ where $H$ is a subgroup of $G$. In doing so, we find there is an equivariant \bok\ spectral sequence which computes the $G$-equivariant homology of $\iota_G^* \thh_H (R)$ and it has on its $E^2$-page the relative theory of Hochschild homology for Green functors given in Definition \ref{defn: relative hh for green functors }. Recall from this definition that the twisted cyclic nerve involved equivariant norms $N^G_H$. Thus once again, we require the use of an $\underline{RO}(G)$-grading to ensure that our grading scheme respects the norm.    

  \begin{thm} \label{thm: relative twisted bok ss}
      Let $H \subseteq G$ be finite subgroups of $S^1$ and let $g=e^{2\pi i/|G|}$ be a generator of $G$. Let $R$ be an $H$-ring spectrum and $E$ a commutative $G$-ring spectrum. Assume that $g$ acts trivially on $E$ and that $\underline{E}_{\ostar}(N^G_H R)$ is flat as a module over $\underline{E}_{\ostar}$. If $R$ has $(\iota_H^*E)$-free homology, then there is a relative twisted \bok\ spectral sequence
      \[
      E_{*,\ostar}^2=\mhh_H^{\underline{E}_{\ostar}, G}(\underline{(\iota_H^*E)}_{\ostar}(R))_* \Rightarrow \underline{E}_{*+\ostar}(\iota^*_G\thh_H(R)).
      \]
  \end{thm}
  
  \begin{proof}
      By Proposition \ref{thm: ekmm simp ss} there is a spectral sequence 
      \[
      E^2_{*, \ostar}= H_*(\underline{E}_{\ostar}(B_{\bullet}^{cy}(N^G_H R; {}^gN^G_H R))) \Rightarrow \underline{E}_{*+\ostar}(|B_{\bullet}^{cy}(N^G_H R; {}^gN^G_H R)|).
      \]
      On the right hand side, this is the twisted cyclic bar construction which defines $\iota^*_G \thh_H(R)$. We wish to identify the $E^2$-page with $\underline{RO}(G)$-graded relative Hochschild homology for Green functors. 
      
We apply the homology functor $\underline{E}_{\ostar}(-)$ level-wise to the relative twisted cyclic bar construction. On $E^1_{p, \ostar}$ this gives
\begin{align*}
     & \underline{E}_{\ostar}({}^gN^G_H R \smsh \overbrace{N^G_H R\smsh ... \smsh N^G_H R}^{p})\\ 
   =\ & \mpi^G_{\ostar}({}^gN^G_H R \smsh N^G_H R \smsh ... \smsh N^G_H R \smsh E) \\ 
  \cong\ & \mpi^G_{\ostar}({}^gN^G_H R \smsh E) \square \mpi^G_{\ostar}(N^G_H R \smsh E) \square ... \square \mpi^G_{\ostar}(N^G_H R \smsh E)\\
   =\ & \underline{E}_{\ostar}({}^gN^G_H R) \square \underline{E}_{\ostar}(N^G_H R) \square ... \square \underline{E}_{\ostar}(N^G_H R)
\end{align*}
where the flatness of $\underline{E}_{\ostar}(N^G_H R)$ as a module over $\underline{E}_{\ostar}$ and an application of the K\"{u}nneth theorem yields the isomorphism. Again, all box products are assumed to be over $\underline{E}_{\ostar}$.  

By the freeness assumption in the statement of the theorem and by Lemma \ref{thm: norm commutes w homology}, we get an isomorphism  
\begin{equation} \label{equn: twisted construction 2}
\underline{E}_{\ostar}(N^G_H R) \cong N_H^G (\underline{(\iota_H^*E)}_{\ostar} R).
\end{equation} 
Combining this isomorphism with the equivalence in Remark \ref{rem: trivial g action} further yields
\begin{equation} \label{eqn: twisted construction 3}
\underline{E}_{\ostar}({}^g N^G_H R) \cong {}^g\underline{E}_{\ostar}(N^G_H R) \cong {}^gN^G_H((\underline{\iota_H^*E})_{\ostar} R). 
\end{equation}
The isomorphisms in \ref{equn: twisted construction 2} and \ref{eqn: twisted construction 3} allow us to conclude that $E^1_{p, \ostar}$ is isomorphic to the following:
\[
{}^gN^G_H((\underline{\iota_H^*E})_{\ostar} R) \square \overbrace{N^G_H((\underline{\iota_H^*E})_{\ostar} R) \square ... \square N^G_H((\underline{\iota_H^*E})_{\ostar} R)}^p.
\]
Term-wise, this is precisely the $p$th level of the twisted cyclic bar construction which defines ${\mhh_H^G}(\underline{(\iota_H^*E)}_{\ostar}(R))$. A diagram chase shows that $d_1$ differential of the spectral sequence agrees with the differential in $\mhh^G_H$, thus we conclude that the relative twisted \bok\ spectral sequence takes the form 
\[
      {E}_{*,\ostar}^2={\mhh_H^{\underline{E}_{\ostar},G}}(\underline{(\iota_H^*E)}_{\ostar}(R))_* \Rightarrow \underline{E}_{*+\ostar}(\iota^*_G\thh_H(R)).
      \]
  \end{proof}


\section{Real Algebraic Structures} \label{chapter: algebraic structures}

Spectral sequence calculations are often quite complex, therefore identifying any additional algebraic structures present in the newly constructed Real \bok\ spectral sequence can provide a computational foothold. In the classical setting, when $A$ a commutative ring spectrum, $\thh(A)$ is an $A$-Hopf algebra in the stable homotopy category (see \cite{mcclure1997thh} and \cite{ekmm}). Angeltveit-Rognes extend this work by constructing simplicial algebraic structure maps on $\thh(A)_{\bullet}$ in \cite{angeltveit2005hopf}. Utilizing the simplicial filtration on $\thh$, the authors then prove that this Hopf algebra structure lifts to the \bok\ spectral sequence. In this section, we use techniques analogous to those of Angeltveit-Rognes to show that $\thr(A)$ is a Hopf algebroid in the $D_2$-equivariant stable homotopy category when $A$ is a commutative ring spectrum with anti-involution. 

\begin{rem}
    If $(A, \omega)$ is a commutative ring spectrum with anti-involution, the distinction of $\omega$ as a map into $A^{op}$ as opposed to $A$ is lost. Hence, we may regard $A$ as simply a commutative $D_2$-ring spectrum with the $D_2$-action of $\omega$. We use this terminology for the remainder of the section.
\end{rem}

The Hopf algebra structure on $\thh$ which is defined in \cite{mcclure1997thh} and \cite{ekmm} is induced by maps on circles since $\thh(R) \cong R \otimes S^1$. Recall from Remark \ref{remark: thr is o2 spectrum} that $\thr(A)$ is an $O(2)$-spectrum. For a nice class of ring spectra with anti-involution, we can recognize $\thr$ as a tensor with $O(2)$.  

\begin{defn} \label{defn: very well pointed}
   An orthogonal $D_2$-spectrum $A$ indexed on a complete universe $\mathcal{U}$ is \textit{well-pointed} if $A(V)$ is well pointed in $\mathit{Top}^{D_2}$ for all finite dimensional orthogonal $D_2$-representations $V$. Further, we say a $D_2$-spectrum $A$ is \textit{very well-pointed} if it is well-pointed and the unit map $S^0 \rightarrow A(\R^0)$ is a Hurewicz cofibration in $\mathit{Top}^{D_2}$.   
\end{defn}

\begin{prop}[\cite{realhh}, Proposition 4.9] \label{thm: thr is tensor with O2}
Let $A$ be a commutative $D_2$-ring spectrum which is very well pointed. Then there is a weak equivalence of $D_2$-spectra
\[
N^{O(2)}_{D_2}A \simeq A \otimes_{D_2} O(2). 
\]
\end{prop}

Note that in the Real equivariant setting, this tensor product occurs over $D_2$. In our case, we utilize the fact that the category of commutative monoids in orthogonal $D_2$-spectra is tensored over the category of $D_2$-sets (see Section 4.1 of \cite{realhh}). More generally, we have that $G$-spectra are tensored in $G$-sets which allows us to define the tensor product of a $G$-spectrum over $G$ as a coequalizer. 

\begin{defn} \label{defn: equivariant tensor product}
  Let $A$ be a commutative $G$-ring spectrum and $X_{\bullet}$ a simplicial $G$-set. The \emph{tensor product over $G$ of $A$ with $X$} is the coequalizer 

\[ \begin{tikzcd}
A \otimes G \otimes X_{\bullet} \arrow[r, "\gamma_1 \otimes id", shift left=1.5ex] \arrow[r, "id \otimes \gamma_2"'] & A \otimes X_{\bullet} \arrow[r] & {A \otimes_{G} X_{\bullet}}
  \end{tikzcd} \]
 where $\gamma_1$ is the $G$-action applied to $A$ and $\gamma_2$ is the $G$-action on $X_{\bullet}$. 
  \end{defn}
  
 The standard simplicial model on $O(2)$ (see $\S 6.3$ of \cite{loday}, for instance) is the geometric realization of a simplicial complex of dihedral groups with the group $D_{2(n+1)}$ at simplicial level $n$ and face and degeneracy maps that commute with the action of $\omega$. To view $\thr(A)$ as the $D_2$-tensor product of $A$ with a simplicial object, we require a further subdivided simplicial model of $O(2)$. In this section, we work with a Segal-Quillen subdivision (see Definition \ref{defn: Segal subdivision}) of this standard simplicial model of $O(2)$. We keep with the conventions of \cite{realhh} and denote this simplicial $D_2$-set by $O(2)_{\bullet}$.
 
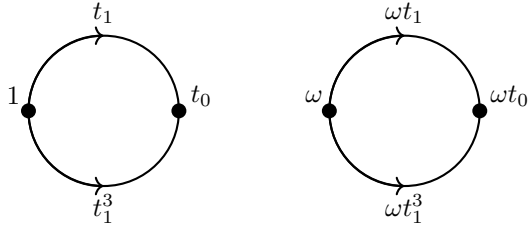
\begin{figure}[h]
    \centering   
\begin{tikzpicture}
    \draw[thick] (0,0) circle(1);
    \node at (-1,0) [circle, fill, inner sep=2pt]{};
    \node at (1,0) [circle, fill, inner sep=2pt]{};
    \node at (-1.2,0.2) {1};
    \node at (1.3,0.2) {$t_0$};
    \draw[->, thick] (-1,0) arc (180:90:1);
    \node at (0,1.3) {$t_1$};
    \draw[->, thick] (-1,0) arc (180:270:1);
    \node at (0,-1.3) {$t_1^3$};

    \draw[thick] (4,0) circle(1);
    \node at (3,0) [circle, fill, inner sep=2pt]{};
    \node at (5,0) [circle, fill, inner sep=2pt]{};
    \node at (2.8,0.2) {$\omega$};
    \node at (5.4,0.2) {$\omega t_0$};
    \draw[->, thick] (3,0) arc (180:90:1);
    \node at (4,1.3) {$\omega t_1$};
    \draw[->, thick] (3,0) arc (180:270:1);
    \node at (4,-1.3) {$\omega t_1^3$};
\end{tikzpicture}
  \caption{The simplicial model $O(2)_{\bullet}$.}
    \label{figure: simplicial O2}
\end{figure}

 One can check that most of the cells in this simplicial object are degenerate, thus upon geometric realization our model of $O(2)_{\bullet}$ can be depicted as two subdivided circles. Explicitly, the simplicial structure of $O(2)_{\bullet}$ at levels 0 and 1 is as follows. At simplicial level 0 it is the group 
    \[
    D_4= \langle t_0, \omega \mid\ t_0^2=1=\omega^2,\ t_0\omega=\omega t_0 \rangle =\{1, t_0, \omega, \omega t_0\}. 
    \]
    To emphasize that $t$ is the generator of the group at the 0th level, we denote it by a subscript $0$. 
    At simplicial level 1 we have 
    \[
    D_8= \langle t_1, \omega \mid\ t_1^4=1=\omega^2,\ t_1 \omega=\omega t_1^3 \rangle =\{1, t_1, t_1^2, t_1^3, \omega, \omega t_1, \omega t_1^2, \omega t_1^3\}.   
    \]
    The elements $1$, $t_1^2$, $\omega$, and $\omega t_1^2$ in simplicial level 1 are in the image of the degeneracy maps. Thus, upon geometric realization to $O(2)_{\bullet}$, we only retain 1-cells indexed by $t_1, t_1^3, \omega t_1,$ and $\omega t_1^3$. This simplicial model $O(2)_{\bullet}$ is depicted in Figure \ref{figure: simplicial O2}. 
  
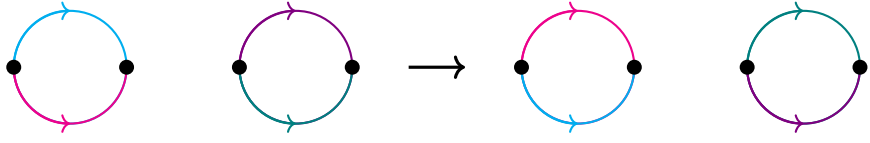
\begin{figure}[h]
    \centering   
\begin{tikzpicture}[scale=0.75]
    \draw[cyan, thick] (0,0) circle(1);
    \draw[->, cyan, thick] (-1,0) arc (180:90:1);
    \draw[->, magenta, thick] (-1,0) arc (180:270:1);
    \draw[magenta, thick] (-1,0) arc (180:358:1);
   \node at (-1,0) [circle, fill, inner sep=2pt]{};
	\node at (1,0) [circle, fill, inner sep=2pt]{};
    \draw[violet, thick] (4,0) circle(1);
    \draw[->, violet, thick] (3,0) arc (180:90:1);
    \draw[->, teal, thick] (3,0) arc (180:270:1);
    \draw[teal, thick] (3,0) arc(180:358:1);
    \node at (3,0) [circle, fill, inner sep=2pt]{};
     \node at (5,0) [circle, fill, inner sep=2pt]{};
     
     \draw[->, very thick] (6,0) -- (7,0);
     
    \draw[magenta, thick] (9,0) circle(1);
    \draw[->, magenta, thick] (8,0) arc (180:90:1);
    \draw[->, cyan, thick] (8,0) arc (180:270:1);
    \draw[cyan, thick] (8,0) arc (180:358:1);
   \node at (8,0) [circle, fill, inner sep=2pt]{};
	\node at (10,0) [circle, fill, inner sep=2pt]{};
    \draw[teal, thick] (13,0) circle(1);
    \draw[->, teal, thick] (12,0) arc (180:90:1);
    \draw[->, violet, thick] (12,0) arc (180:270:1);
    \draw[violet, thick] (12,0) arc(180:358:1);
    \node at (12,0) [circle, fill, inner sep=2pt]{};
     \node at (14,0) [circle, fill, inner sep=2pt]{};

\end{tikzpicture}
  \caption{The action of $t$ on $O(2)_{\bullet}$.}
    \label{figure: t action}
\end{figure}

To simplify the remaining figures in this section, we will cease to label every cell. This subdivided model of $O(2)_{\bullet}$ is endowed with a $D_4$-action where $\omega$ acts by swapping the two circles and $t$ reflects within each circle (see Figure \ref{figure: t action}). In the equivariant tensor product $A \otimes_{D_2} O(2)_{\bullet}$, the tensor occurs over the action of $t$. 

  We now use this simplicial structure to define algebraic structure maps on Real topological Hochschild homology.

\begin{lemma} \label{thm: thr is an algebra}
    Let $A$ be a commutative $D_2$-ring spectrum. The Real topological Hochschild homology of $A$ is a commutative $A$-algebra in $D_2$-spectra. 
\end{lemma}

\begin{proof}

  The unit map $\eta: A \rightarrow \thr(A)$ is induced by the inclusion of ${D_2}_{\bullet}$ into $O(2)_{\bullet}$ as the zero cells $1$ and $\omega$. In order to define a product map, we must first recognize $\thr(A) \smsh_A \thr(A)$ as a $D_2$-tensor with a simplicial set. Consider the following pushout in $D_2$-simplicial sets:
\[\begin{tikzcd}
	 {D_2}_{\bullet} \arrow[r] \arrow[d]& {O(2)_{\bullet}} \arrow[d]\\
	 O(2)_{\bullet} \arrow[r] & O(2)_{\bullet} \vee_{D_2} O(2)_{\bullet}.	
\end{tikzcd}\]
The simplicial object $O(2)_{\bullet} \vee_{D_2} O(2)_{\bullet}$ in the bottom right corner of this diagram is depicted in Figure \ref{figure: o2 wedge o2}. 

\begin{figure}[h]
    \centering   
\begin{tikzpicture}
    \draw[thick] (0,0) circle(1);
    \node at (-1,0) [circle, fill, inner sep=2pt]{};
    \node at (1,0) [circle, fill, inner sep=2pt]{};
    \draw[->, thick] (1,0) arc (0:90:1);
    \draw[->, thick] (1,0) arc (0:-90:1);
     \draw[thick] (2,0) circle(1);
     \node at (3,0) [circle, fill, inner sep=2pt]{};
      \draw[->, thick] (1,0) arc (180:90:1);
    \draw[->, thick] (1,0) arc (180:270:1);

    \draw[thick] (5,0) circle(1);
    \node at (4,0) [circle, fill, inner sep=2pt]{};
    \node at (6,0) [circle, fill, inner sep=2pt]{};
    \draw[->, thick] (6,0) arc (0:90:1);
    \draw[->, thick] (6,0) arc (0:-90:1);
     \draw[thick] (7,0) circle(1);
     \node at (8,0) [circle, fill, inner sep=2pt]{};
      \draw[->, thick] (6,0) arc (180:90:1);
    \draw[->, thick] (6,0) arc (180:270:1);
\end{tikzpicture}
  \caption{The simplicial object $O(2)_{\bullet} \vee_{D_2} O(2)_{\bullet}$.}
    \label{figure: o2 wedge o2}
\end{figure}
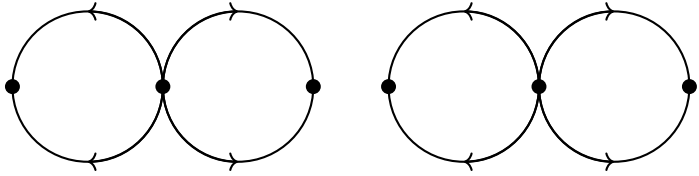

 The tensor product in $D_2$-spectra preserves pushouts so upon applying the functor $|A \otimes_{D_2} (-)|$ to this diagram we obtain
  \[\begin{tikzcd}
	 A \arrow[r] \arrow[d]& \thr(A) \arrow[d]\\
	 \thr(A) \arrow[r] & {|A \otimes_{D_2}(O(2)_{\bullet} \vee_{D_2} O(2)_{\bullet})|.}	
\end{tikzcd}\]
 The pushout of this diagram defines the relative smash product so we have identified $|A \otimes_{D_2} (O(2)_{\bullet} \vee_{D_2} O(2)_{\bullet})|$ as  $\thr(A) \smsh_A \thr(A)$. 
  
  Now we may define a product map
    \[
    \mu: \thr(A) \smsh_A \thr(A) \rightarrow \thr(A)
    \]
    which is induced by the wedge of two copies of $O(2)_{\bullet}$ that folds one copy of $O(2)_{\bullet}$ onto the other, as in Figure \ref{figure: thr product}. For clarity, the second wedge of circles is not color coded in this figure but the identifications are made in the same way.  
    
     \begin{figure}[hbt!]
    \centering   
\begin{tikzpicture}[scale=.65]
    \draw[thick, red] (-1,0) arc (180:0:1);
    \draw[thick, blue] (-1,0) arc(180:360:1);
    \draw[thick, olive] (1,0) arc (180:0:1);
    \draw[thick, orange] (1,0) arc (180: 360:1);
    \node at (-1,0) [circle, fill, cyan, inner sep=2pt]{};
    \node at (1,0) [circle, fill, violet, inner sep=2pt]{};
    \node at (1,0) [circle, fill, teal, inner sep=1pt]{};
    \node at (3,0) [circle, fill, magenta, inner sep=2pt]{};

    \node at (3,-1.5) {$O(2)_{\bullet} \vee_{D_2} O(2)_{\bullet}$};

    \draw[thick] (5,0) circle(1);
    \draw[thick] (7,0) circle(1);

    \draw[->, thick, dashed] (9,0) -- (11,0);
    \node at (10,0.3) {$\mu$};

    \draw[thick] (16,0) circle(1);

    \draw[thick, red] (12,0) arc (180:0:1);
    \draw[thick, olive] (12.05,0) arc (180:0:.95);
    \draw[thick, blue] (12,0) arc (180:360:1);
    \draw[thick, orange] (12.05, 0) arc (180:360:.95);
    \node at (12,0) [circle, fill, cyan, inner sep=2pt]{};
     \node at (12,0) [circle, fill, magenta, inner sep=1pt]{};
    \node at (14,0) [circle, fill, violet, inner sep=2pt]{};
    \node at (14,0) [circle, fill, teal,  inner sep=1pt]{};

    \node at (15,-1.5) {$O(2)_{\bullet}$};
\end{tikzpicture}
  \caption{The simplicial map $\mu$ which induces a product.}
    \label{figure: thr product}
\end{figure}
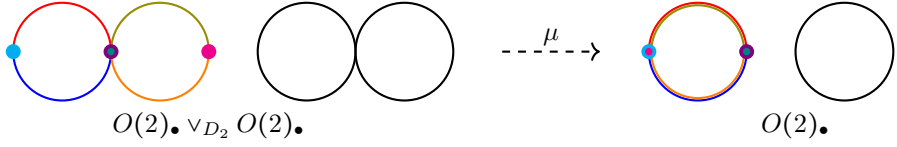
    
 Note that both the product and unit are equivariant with respect to the actions of $\omega$ and $t$ on $O(2)_{\bullet}$. We omit the details, but one can check that these maps satisfy commutative diagrams in $D_2$-simplicial sets in order to show that the $A$-algebra structure on $\thr(A)$ is unital and associative.
  \end{proof} 
  
  In the preceding proof we defined a unit map on $\simpotwo$ by including $D_2$ as the cells $1$ and $\omega$ of $O(2)_{\bullet}$. However, our simplicial model of $O(2)$ also includes another pair of points that one could have used instead. The presence of two distinct unit maps suggests that rather than inheriting a Hopf algebra structure like that of $\thh$, Real topological Hochschild homology has the structure of a more general object called a Hopf algebroid. 
  \begin{defn} \label{defn: hopf algebroid}
    A \emph{Hopf algebroid} over a commutative ring $k$ is a pair of commutative $k$-algebras $(A, R)$ together with:
    \begin{itemize}
        \item a left unit map $\eta_L: A \rightarrow R$
        \item a right unit map $\eta_R: A \rightarrow R$
        \item a coproduct map $\delta: R \rightarrow R \otimes_A R$
        \item a counit map $\varepsilon: R \rightarrow A$
        \item an antipode map $\chi: R \rightarrow R$ which squares to the identity
    \end{itemize}
    satisfying the following conditions:
    
    \begin{enumerate}
    \item $\varepsilon \circ \eta_L = \varepsilon \circ \eta_R = id_A$
    \item $\varepsilon$ is counital
    \item $\delta$ is coassociative
    \item $\chi \circ \eta_R = \eta_L$ and $\chi \circ \eta_L = \eta_R$
    \item There exist maps $\mu_R$ and $\mu_L$ such that the following diagram commutes, where $\varphi$ is the multiplication map on $R \otimes R$ as a $k$-algebra:
    \[\begin{tikzcd}
    A \arrow[d, "\eta_R"] & R \arrow[l, "\varepsilon"'] \arrow[r, "\varepsilon"] \arrow[d, "\delta"] & A \arrow[d, "\eta_L"] \\
    R & R \otimes_A R \arrow[l, dashed, "\mu_R"'] \arrow[r, dashed, "\mu_L"] & R \\
    R \otimes_k R \arrow[u, "\varphi"'] & R \otimes_k R \arrow[l, "\chi \otimes id"'] \arrow[r, "id \otimes \chi"] \arrow[u] & R \otimes_k R. \arrow[u, "\varphi"] 
        \end{tikzcd}\]
    \end{enumerate}
    \end{defn}
    
    Hopf algebroids generalize the notion of a Hopf algebra akin to the generalization of a group via a groupoid. When the left and right units coincide, $R$ is simply an $A$-Hopf algebra. 
    
  We claim that $\thr(A)$ has a Hopf algebroid structure in the $D_2$-homotopy category. To show the existence of this structure we will define a coproduct map on $\simpotwo$. The most natural map is a pinch map, however attempting to identify the zero cells $1 \sim t_0$ and $\omega \sim \omega t_0$ in Figure \ref{figure: simplicial O2} will not produce $\thr(A) \smsh_A \thr(A)$ upon tensoring with $A$ over $D_2$. Angeltveit-Rognes encounter a similar problem in trying to define a map of simplicial circles that lifts to a coproduct on $\thh$. Their solution is to work with a subdivided model of the simplicial circle and define a coproduct up to homotopy. We adopt this method in order to define a coproduct on $\thr$, beginning with an additional subdivision of $O(2)_{\bullet}$.    
  
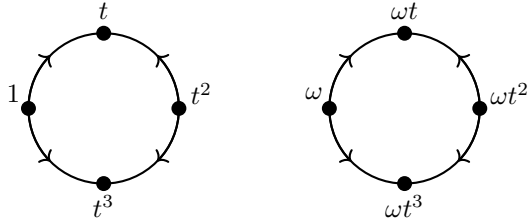
\begin{figure}[h]
    \centering   
\begin{tikzpicture}
    \draw[thick] (0,0) circle(1);
    \node at (-1,0) [circle, fill, inner sep=2pt]{};
    \node at (1,0) [circle, fill, inner sep=2pt]{};
    \node at (-1.2,0.2) {1};
    \node at (1.3,0.2) {$t^2$};
    \node at (0,1) [circle, fill, inner sep=2pt]{};
    \node at (0,-1) [circle, fill, inner sep=2pt]{};
    \draw[->, thick] (-1,0) arc (180:135:1);
    \node at (0,1.3) {$t$};
    \draw[->, thick] (-1,0) arc (180:225:1);
    \node at (0,-1.3) {$t^3$};
    \draw[->, thick] (1,0) arc (0:45:1);
	\draw[->, thick] (1,0) arc (0:-45:1);
    \draw[thick] (4,0) circle(1);
    \node at (3,0) [circle, fill, inner sep=2pt]{};
    \node at (5,0) [circle, fill, inner sep=2pt]{};
    \node at (4,1) [circle, fill, inner sep=2pt]{};
    \node at (4,-1) [circle, fill, inner sep=2pt]{};
    \node at (2.8,0.2) {$\omega$};
    \node at (5.4,0.2) {$\omega t^2$};
    \draw[->, thick] (3,0) arc (180:225:1);
    \draw[->, thick] (3,0) arc (180:135:1);
    \node at (4,1.3) {$\omega t$};
    \draw[->, thick] (5,0) arc (0:45:1);
     \draw[->, thick] (5,0) arc (0:-45:1);
    \node at (4,-1.3) {$\omega t^3$};
\end{tikzpicture}
 \caption{The double model $dO(2)_{\bullet}$.}
  \label{figure: double thr}
\end{figure}

  \begin{defn}
  The \emph{double model} of $O(2)_{\bullet}$ is the Segal-Quillen subdivision
\[
dO(2)_{\bullet} := \text{sq} (O(2)_{\bullet}) = \text{sq} (\text{sq} (D_{2(\bullet+1)})), 
\]
depicted geometrically in Figure \ref{figure: double thr}. We refer to the tensor product $A \otimes_{D_2} dO(2)_{\bullet}$ as the double model of $\thr(A)$. 

\end{defn}
  
We must now verify that $d\thr(A) \simeq \thr(A)$. To do so, we argue that $\thr$ and $d\thr$ can be understood as pushouts by taking the $D_2$-tensor of diagrams in simplicial sets and then demonstrate an equivalence of pushouts. We begin by recognizing $\thr(A)$ as a pushout.

For a ring spectrum $A$, we let $B(A)_{\bullet}$ denote the double bar construction $B(A, A, A)_{\bullet}$. For the sake of clarity, we label copies of the ring spectrum $A$ as $A_i$. A Segal-Quillen subdivision of $B(A)_{\bullet}$, denoted by $\text{sq}B(A)_{\bullet}$, doubles the number of copies of $A$ at each level:
\[\text{sq} B(A)_n = A_0 \smsh A_1 \smsh ... \smsh A_{2n} \smsh A_{2n+1} .\] Subdividing also induces $D_2$-action on $\text{sq} B(A)_n$ which swaps $A_i$ and $A_{2n+1-i}$. In particular, note that the coefficients $A_0$ and $A_{2n+1}$ are exchanged by the $D_2$-action on $\text{sq} B(A)_{\bullet}$.  

We define $\text{sq}B(A)$ to be the geometric realization $ | \text{sq} B(A)_{\bullet} |$ and similarly write $B(A) = | B(A, A, A)_{\bullet} |$. 

 \begin{prop} \label{thm: thr is a pushout}
    Let $A$ be a commutative $D_2$-ring spectrum. Then $\thr(A)$ is the pushout in $D_2$-spectra of the diagram
     \begin{equation} \label{diagram: thr pushout}
  \begin{tikzcd}
 \text{sq} B(A) & A \smsh A \arrow[l] \arrow[r] & A.
  \end{tikzcd}
  \end{equation}
 where the left map includes $A \smsh A$  as the coefficients in the subdivided bar construction and the map on the right is multiplication.
\end{prop}

\begin{proof}
First, we consider a pushout diagram in $D_2$-simplicial sets. Let $\Delta^1$ denote the standard 1-simplex. In the diagram, 
\begin{equation}\label{diagram: simplicial thr pushout}
\begin{tikzcd}
{D_4}_{\bullet} \arrow[r] \arrow[d] & {D_2}_{\bullet} \arrow[d] \\
D_2 \otimes \text{sq} \Delta^1_{\bullet} \arrow[r] & O(2)_{\bullet},
\end{tikzcd}
\end{equation}
the top map indentifies $1$ and $t$ in $D_4$ and the map on the left includes $D_4$ as the boundaries of the two subdivided 1-simplices. 

We apply the functor $|A \otimes_{D_2} (-) |$ to this diagram. Since $| A \otimes_{D_2} {D_4}_{\bullet} |$ is $A \smsh A$ endowed with the swap action, all that remains is to identify $| A \otimes_{D_2} (D_2 \otimes \text{sq} \Delta^1_{\bullet})|$ with $\text{sq} B(A)$. At simplicial level $n$, $\text{sq} \Delta^1$ is the set $\Delta^1_{2n+1}=\{x_0, x_1, ..., x_{2(n+1)}\}$ with a $D_2$-action given swapping $x_i$ and $x_{2(n+1)-i}$. A level-wise comparison shows that $A \otimes (\text{sq} \Delta^1)_n= (\text{sq} B(A))_n$ and one may check that the face and degeneracy maps agree. 

Tensoring with $A$ over $D_2$ preserves pushouts, hence $\thr(A) \simeq | A \otimes_{D_2} O(2)_{\bullet}|$ arises as a pushout from the diagram in the statement of the proposition. 
\end{proof}

We now employ a similar technique to recognize the double model of $\thr$ as arising from a simplicial pushout. 

\begin{prop} \label{thm: dthr is a pushout}
    Let $A$ be a commutative $D_2$-ring spectrum. Then $d\thr(A)$ is the pushout in $D_2$-spectra of the diagram
    \begin{equation} \label{diagram: dthr pushout}
  \begin{tikzcd}
 \text{sq} B(A) & A \smsh A \arrow[l] \arrow[r] & \text{sq} B(A).
  \end{tikzcd}
  \end{equation}
 where both maps include $A \smsh A$  as the coefficients in the subdivided bar construction.
 \end{prop}
 
 \begin{proof}
 In this case, the correct pushout of $D_2$-simplicial sets is the following diagram:
 \begin{equation}\label{diagram: simplicial dthr pushout}
\begin{tikzcd}
{D_4}_{\bullet} \arrow[r] \arrow[d] & D_2 \otimes \text{sq} \Delta^1_{\bullet} \arrow[d] \\
D_2 \otimes \text{sq} \Delta^1_{\bullet} \arrow[r] & dO(2)_{\bullet}.
\end{tikzcd}
\end{equation}
We apply the functor $| A \otimes_{D_2} (-)|$ to the entire diagram and make the same identifications as in the preceding proof. Pushouts are preserved under the tensor hence we obtain the pushout in the statement of the proposition. 
\end{proof}

\begin{figure}[h]
    \centering   
\begin{tikzpicture}[scale=.75]
    \draw[thick] (0,0) circle(1);
    \node at (-1,0) [circle, fill, inner sep=2pt, color=teal]{};
    \node at (0,1) [circle, fill, inner sep=2pt, color=cyan]{};
    \node at (0,-1) [circle, fill, inner sep=2pt, color=cyan]{};
    \node at (1,0) [circle, fill, inner sep=2pt, color=cyan]{};
    \draw[thick,cyan] (1,0) arc (0:90:1);
    \draw[thick,cyan] (1,0) arc (0:-90:1);
    \node at (-1.4,0) {1};

    \draw[thick] (3,0) circle(1);
    \node at (2,0) [circle, fill, inner sep=2pt, color=violet]{};
    \node at (3,1) [circle, fill, inner sep=2pt, color=magenta]{};
    \node at (3,-1) [circle, fill, inner sep=2pt, color=magenta]{};
    \node at (4,0) [circle, fill, inner sep=2pt, color=magenta]{};
    \draw[thick, magenta] (4,0) arc (0:90:1);
    \draw[thick, magenta] (4,0) arc (0:-90:1);
    \node at (1.6,0) {$\omega$};

    \draw[thick, ->] (4.5,0) -- (6.5,0);
    \node at (5.5,0.3) {$\pi$};

    \draw[thick] (8,0) circle(1);
    \node at (7,0) [circle, fill, inner sep=2pt, color=teal]{};
    \node at (9,0) [circle, fill, inner sep=2pt, color=cyan]{};

    \draw[thick] (11,0) circle(1);
    \node at (10,0) [circle, fill, inner sep=2pt, color=violet]{};
    \node at (12,0) [circle, fill, inner sep=2pt, color=magenta]{};
\end{tikzpicture}
 \caption{The simplicial homotopy $\pi$.}
  \label{figure: thr homotopy}
\end{figure}
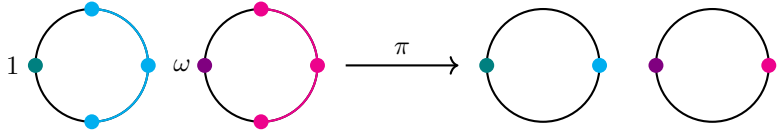

\begin{lemma}\label{thm: double thr is htpy equiv to thr}
    Let $A$ be a commutative $D_2$-ring spectrum which is cofibrant as a $D_2$-spectrum. Then there is a $D_2$-weak equivalence 
  \[
  \pi: d\thr(A) \xrightarrow{\simeq} \thr(A)
  \]
 which is induced by the simplicial homotopy collapsing one half of each circle in $\simpotwo$ to a point, depicted in Figure \ref{figure: thr homotopy}. 
\end{lemma}

 \begin{proof}
 We can construct the following commutative diagram in $D_2$-spectra from the pushouts in Proposition \ref{thm: thr is a pushout} and Proposition \ref{thm: dthr is a pushout}: 
 \begin{equation} \label{eqn: thr dthr homotopy diagram}
\begin{tikzcd}
\text{sq}B(A) \arrow[d, equal] & A \smsh A \arrow[l] \arrow[r] \arrow[d, equal] & \text{sq} B(A) \arrow[d, "\simeq"]\\
\text{sq}B(A) & A \smsh A \arrow[l] \arrow[r] & A. 
\end{tikzcd}
\end{equation}

The homotopy equivalence $B(A, A, A) \simeq A$ (see \cite{ekmm}, IV. 7.3 and XII.1.2) arises simplicially in the bar construction by treating $A$ as a constant simplicial object, $A_{\bullet}$. Since the bar construction we consider here is a subdivision of $B(A)_{\bullet}$, the map inducing the equivalence $B(A)_{\bullet} \rightarrow A_{\bullet}$ at level $2n+1$ induces the equivalence $\text{sq} B(A)_{\bullet} \rightarrow A_{\bullet}$ at level $n$. The homotopy is given by an iterated composites of unit and multiplication maps, which are all $D_2$-equivariant maps. Since pushouts preserve weak equivalences by \cite{ekmm} III.8.2, we obtain a weak equivalence between the pushout along the top row and the pushout along the bottom, hence $d\thr(A) \xrightarrow{\simeq} \thr(A)$. 

We verify that this is a $D_2$-weak equivalence by checking that the map is a weak equivalence on the geometric fixed points. The geometric fixed points functor, denoted by $\Phi^{D_2}$, commutes with colimits and since $\thr(A)$ and $d\thr(A)$ are both colimits by Propositions \ref{thm: thr is a pushout} and \ref{thm: dthr is a pushout}, we have that
\begin{align*}
    \Phi^{D_2}(d\thr(A))&\cong \colim(\Phi^{D_2} (\text{sq} B(A)) \leftarrow \Phi^{D_2}(A \smsh A) \rightarrow \Phi^{D_2} (\text{sq} B(A)))\\
     \Phi^{D_2}(\thr(A))&\cong \colim(\Phi^{D_2} (\text{sq} B(A)) \leftarrow \Phi^{D_2}(A \smsh A) \rightarrow \Phi^{D_2} (A)).
\end{align*}
To compare the terms on the right, we recall that $\Phi^{D_2}$ commutes with the smash product and that this functor is applied level-wise to a simplicial object. Hence $\Phi^{D_2}(\text{sq} B(A)) \cong \text{sq} B(\Phi^{D_2} A)$. Since a two-sided bar construction and its right coefficients are weakly equivalent we have that $\text{sq} B(\Phi^{D_2} A) \simeq \Phi^{D_2} A$. Thus we conclude that the equivalence $d\thr(A) \simeq A$ is a $D_2$-weak equivalence of spectra. 

\end{proof}

Equipped with this homotopy equivalence between our two models of $\thr$, we are now able to describe a coproduct structure and thus show that $\thr$ has the structure of a Hopf algebroid.  

\begin{thm} \label{thm: thr is a hopf algebroid}
    For a commutative $D_2$-ring spectrum $A$, $\thr(A)$ is a Hopf algebroid in the $D_2$-equivariant stable homotopy category. 
\end{thm}

\begin{proof}
The left unit map $\eta_L: A \rightarrow \thr(A)$ is induced by including $D_2$ as the zero cells $1$ and $\omega$ in $\simpotwo$ while the right $\eta_R$ includes $D_2$ as $t$ and $\omega t$. The counit $\varepsilon: \thr(A) \rightarrow A$ collapses the first circle in $\simpotwo$ to $1$ and the second circle to $\omega$. 

As noted, we must work with $d\thr(A)$ in order to define a coproduct map. Consider the map $\delta': d\thr(A) \rightarrow \thr(A) \smsh_A \thr(A)$ induced by the simplicial map on $d\simpotwo$ which identifies $t \sim t^3$ and $\omega t \sim \omega t^3$. This is depicted in Figure \ref{figure: thr coproduct}. 

 \begin{figure}[h]
    \centering 
\begin{tikzpicture}[scale=.7]
    \draw[thick] (0,0) circle(1);
   
    \draw[->, thick] (-1,0) arc (180:135:1);
    \draw[->, thick] (-1,0) arc (180:225:1);
    \draw[->, thick] (1,0) arc (0:-45:1);
    \draw[->, thick] (1,0) arc (0:45:1);
     \node at (-1,0) [circle, fill, inner sep=2pt, color=blue]{};
    \node at (0,1) [circle, fill, inner sep=2pt, color=magenta]{};
    \node at (0,-1) [circle, fill, inner sep=2pt, color=cyan]{};
    \node at (1,0) [circle, fill, inner sep=2pt, color=red]{};
    \node at (-1.4,0) {1};

    \draw[thick] (3,0) circle(1);
    
    \draw[->, thick] (2,0) arc (180:135:1);
    \draw[->, thick] (2,0) arc (180:225:1);
    \draw[->, thick] (4,0) arc (0:-45:1);
    \draw[->, thick] (4,0) arc (0:45:1);
    \node at (2,0) [circle, fill, inner sep=2pt, color=yellow]{};
    \node at (3,1) [circle, fill, inner sep=2pt, color=teal]{};
    \node at (3,-1) [circle, fill, inner sep=2pt, color=violet]{};
    \node at (4,0) [circle, fill, inner sep=2pt, color=orange]{};
    \node at (1.6,0) {$\omega$};

    \node at (2,-1.6) {$O(2)_{\bullet}$};

    \draw[thick, ->] (4.5,0) -- (6.5,0);
    \node at (5.5,0.3) {$\delta'$};

    \draw[thick] (8,0) circle(1);
    \draw[thick] (10,0) circle(1);
    \draw[->,thick] (7,0) arc (180:90:1);
    \draw[->,thick] (7,0) arc (180:270:1);
    \draw[->, thick] (11,0) arc (0:90:1);
    \draw[->, thick] (11,0) arc (0:-90:1);
    \node at (7,0) [circle, fill, inner sep=2pt, color=blue]{};
    \node at (9,0) [circle, fill, inner sep=2pt, color=cyan]{};
    \node at (9,0) [circle, fill, inner sep=1pt, color=magenta]{};
    \node at (11,0) [circle, fill, inner sep=2pt, color=red]{};

    \node at (11.6,-1.6) {$O(2)_{\bullet} \lor_{D_2} O(2)_{\bullet}$};

    \draw[thick] (13,0) circle(1);
    \draw[thick] (15,0) circle(1);
    \draw[->,thick] (12,0) arc (180:90:1);
    \draw[->,thick] (12,0) arc (180:270:1);
    \draw[->, thick] (16,0) arc (0:90:1);
    \draw[->, thick] (16,0) arc (0:-90:1);
    \node at (12,0) [circle, fill, inner sep=2pt, color=yellow]{};
    \node at (14,0) [circle, fill, inner sep=2pt, color=violet]{};
    \node at (14,0) [circle, fill, inner sep=1pt, color=teal]{};
    \node at (16,0) [circle, fill, inner sep=2pt, color=orange]{};


\end{tikzpicture}
 \caption{The simplicial map inducing a coproduct.}
  \label{figure: thr coproduct}
\end{figure}
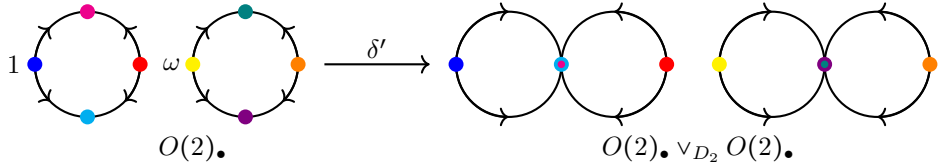

The coproduct on $\thr$ is then given by the composite 
\[
\delta: \thr(A) \xrightarrow{\pi^{-1}} d\thr(A) \xrightarrow{\delta'} \thr(A) \smsh_A \thr(A), 
\]
where $\pi$ is the equivalence from Lemma \ref{thm: double thr is htpy equiv to thr}.

Finally, we define an antipodal map on $\thr$ induced by the map on $O(2)_{\bullet}$ which reflects each circle across a vertical axis. The action of this antipodal map, which we will call $\chi$, is depicted in Figure \ref{figure: thr antipode}.

 \begin{figure}[h]
    \centering   
\begin{tikzpicture}[scale=0.75]
    \draw[thick] (0,0) circle(1);
    \draw[->,thick] (-1,0) arc (180:90:1);
    \draw[->, thick] (-1,0) arc (180:270:1);
   \node at (-1,0) [circle, fill, cyan, inner sep=2pt]{};
	\node at (1,0) [circle, fill, magenta, inner sep=2pt]{};
    \draw[thick] (4,0) circle(1);
    \draw[->,thick] (3,0) arc (180:90:1);
    \draw[->,thick] (3,0) arc (180:270:1);
    \node at (3,0) [circle, fill, violet, inner sep=2pt]{};
     \node at (5,0) [circle, fill, teal, inner sep=2pt]{};
     
     \draw[->, very thick] (6,0) -- (7,0);
     \node at (6.5, 0.4) {$\chi$};
          
    \draw[thick] (9,0) circle(1);
    \draw[->,thick] (10,0) arc (0:90:1);
    \draw[->, thick] (10,0) arc (0:-90:1);
   \node at (8,0) [circle, fill, magenta, inner sep=2pt]{};
	\node at (10,0) [circle, fill, cyan, inner sep=2pt]{};
    \draw[thick] (13,0) circle(1);
    \draw[->, thick] (14,0) arc (0:90:1);
    \draw[->, thick] (14,0) arc (0:-90:1);
    \node at (12,0) [circle, fill, teal, inner sep=2pt]{};
     \node at (14,0) [circle, fill, violet, inner sep=2pt]{};

\end{tikzpicture}
  \caption{The action of the map $\chi$ which induces an antipode on $\thr$.}
    \label{figure: thr antipode}
\end{figure}
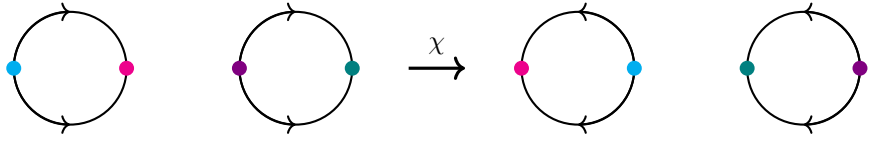

To verify that these structure maps satisfy the commutativity relations of a Hopf algebroid as stated in Definition \ref{defn: hopf algebroid}, we check that the relations hold for the simplicial maps of $O(2)_{\bullet}$. From a visual inspection it is clear that the antipodal map swaps the units. That composing $\varepsilon$ with either of the units produces the identity map is also clear.  

To verify counitality and coassociativity, we may check that the appropriate diagrams commute in $D_2$-simplicial sets. For example, the following diagram demonstrates that $\varepsilon$ is counital,
\[\begin{tikzcd}
d\simpotwo \arrow[r, "\delta'"] \arrow[d, "\delta'"] \arrow[ddr, "\simeq"]  & \simpotwo \vee_{D_2} \simpotwo \arrow[dd,"id \vee \varepsilon"] \\
\simpotwo \vee_{D_2} \simpotwo \arrow[d, "\varepsilon \vee id"] & \\
\simpotwo \arrow[r, "\simeq"] & \simpotwo
\end{tikzcd}\]
where the composition on the left is a simplicial homotopy that collapses $t, t^3$ to $1$ and $\omega t, \omega t^3$ to $\omega$ in Figure \ref{figure: double thr}. Note, this map collapses the other half of the circles as opposed to $\pi$ depicted in Figure \ref{figure: thr homotopy}. This homotopy is represented by the diagonal arrow. We omit the diagram, but one may similarly check that the coassociativity relation holds.    

Finally, we show the existence of two maps (denoted by $\mu_L$ and $\mu_R$ so as to suggest a left and right multiplication) which make the following diagram commute:  

\begin{equation} \label{diagram: thr compatibility}
\begin{tikzcd}
 {D_2}_{\bullet} \arrow[d, "\eta_R"] & d\simpotwo \arrow[l, "\varepsilon \circ \pi"'] \arrow[d, "\delta'"] \arrow[r, "\varepsilon \circ \pi"] & {D_2}_{\bullet} \arrow[d, "\eta_L"] \\
 \simpotwo & \simpotwo \vee_{D_2} \simpotwo \arrow[l, dashed, "\mu_R"'] \arrow[r, dashed, "\mu_L"] & \simpotwo \\
 \simpotwo \sqcup \simpotwo \arrow[u, "\varphi"] &  \simpotwo \sqcup \simpotwo \arrow[u] \arrow[l, "\chi \sqcup id"'] \arrow[r, "id \sqcup \chi"] & \simpotwo \sqcup \simpotwo. \arrow[u, "\varphi"]
 \end{tikzcd}
 \end{equation}

Here, the map $\varphi$ associates the corresponding parts of two copies of $\simpotwo$. The maps $\mu_L$ and $\mu_R$ both fold the wedged copies of $\simpotwo$ onto themselves but in different directions; $\mu_L$ folds to the left and $\mu_R$ folds to the right. Both maps are depicted in Figure \ref{figure: thr left product} and Figure \ref{figure: thr right product} with overlapping colors denoting cells that get associated under the maps. In these pictures, we do not color code the components of the second circle in $\simpotwo$ but the associations are made in the same way as depicted in the first circle. 

 \begin{figure}[hbt!]
    \centering   
\begin{tikzpicture}[scale=.65]
    \draw[thick, red] (-1,0) arc (180:0:1);
    \draw[thick, blue] (-1,0) arc(180:360:1);
    \draw[thick, olive] (1,0) arc (180:0:1);
    \draw[thick, orange] (1,0) arc (180: 360:1);
    \node at (-1,0) [circle, fill, cyan, inner sep=2pt]{};
    \node at (1,0) [circle, fill, violet, inner sep=2pt]{};
    \node at (1,0) [circle, fill, teal, inner sep=1pt]{};
    \node at (3,0) [circle, fill, magenta, inner sep=2pt]{};

    \node at (3,-1.5) {$O(2)_{\bullet} \vee_{D_2} O(2)_{\bullet}$};

    \draw[thick] (5,0) circle(1);
    \draw[thick] (7,0) circle(1);

    \draw[->, thick, dashed] (9,0) -- (11,0);
    \node at (10,0.3) {$\mu_L$};

    \draw[thick] (16,0) circle(1);

    \draw[thick, red] (12,0) arc (180:0:1);
    \draw[thick, olive] (12.05,0) arc (180:0:.95);
    \draw[thick, blue] (12,0) arc (180:360:1);
    \draw[thick, orange] (12.05, 0) arc (180:360:.95);
    \node at (12,0) [circle, fill, cyan, inner sep=2pt]{};
     \node at (12,0) [circle, fill, magenta, inner sep=1pt]{};
    \node at (14,0) [circle, fill, violet, inner sep=2pt]{};
    \node at (14,0) [circle, fill, teal,  inner sep=1pt]{};

    \node at (15,-1.5) {$O(2)_{\bullet}$};
\end{tikzpicture}
  \caption{The simplicial map $\mu_L$ which acts as a left-product.}
    \label{figure: thr left product}
\end{figure}
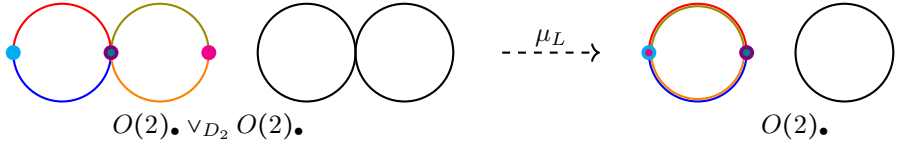

 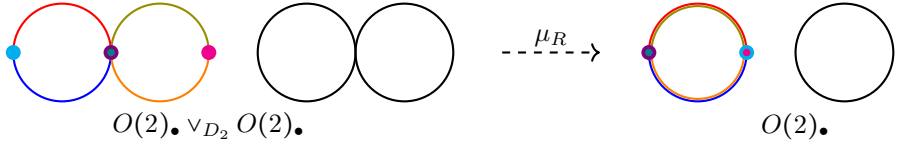
\begin{figure}[h]
    \centering   
\begin{tikzpicture}[scale=.65]
    \draw[thick, red] (-1,0) arc (180:0:1);
    \draw[thick, blue] (-1,0) arc(180:360:1);
    \draw[thick, olive] (1,0) arc (180:0:1);
    \draw[thick, orange] (1,0) arc (180: 360:1);
    \node at (-1,0) [circle, fill, cyan, inner sep=2pt]{};
    \node at (1,0) [circle, fill, violet, inner sep=2pt]{};
    \node at (1,0) [circle, fill, teal, inner sep=1pt]{};
    \node at (3,0) [circle, fill, magenta, inner sep=2pt]{};

    \node at (3,-1.5) {$O(2)_{\bullet} \vee_{D_2} O(2)_{\bullet}$};

    \draw[thick] (5,0) circle(1);
    \draw[thick] (7,0) circle(1);

    \draw[->, thick, dashed] (9,0) -- (11,0);
    \node at (10,0.3) {$\mu_R$};

    \draw[thick] (16,0) circle(1);

    \draw[thick, red] (12,0) arc (180:0:1);
    \draw[thick, olive] (12.05,0) arc (180:0:.95);
    \draw[thick, blue] (12,0) arc (180:360:1);
    \draw[thick, orange] (12.05, 0) arc (180:360:.95);
    \node at (12,0) [circle, fill, violet, inner sep=2pt]{};
     \node at (12,0) [circle, fill, teal, inner sep=1pt]{};
    \node at (14,0) [circle, fill, cyan, inner sep=2pt]{};
    \node at (14,0) [circle, fill, magenta,  inner sep=1pt]{};

    \node at (15,-1.5) {$O(2)_{\bullet}$};
\end{tikzpicture}
  \caption{The simplicial map $\mu_R$ which acts as a right-product.}
    \label{figure: thr right product}
\end{figure}

We claim composites $\mu_L \circ \delta '$ and  $\mu_R \circ \delta'$ in \ref{diagram: thr compatibility} both factor through the map
\[
\partial \bar{\Delta}^2 \sqcup \partial \bar{\Delta}^2  \hookrightarrow  \bar{\Delta}^2 \sqcup \bar{\Delta}^2 \rightarrow \bar{\Delta}^1 \sqcup \bar{\Delta}^1 \rightarrow \simpotwo,
\]
where $\bar{\Delta}^2$ denotes the subdivided Real 2-simplex with a $D_2$-action that reflects across the vertical axis. This composition is depicted in Figure \ref{figure: thr contraction}.

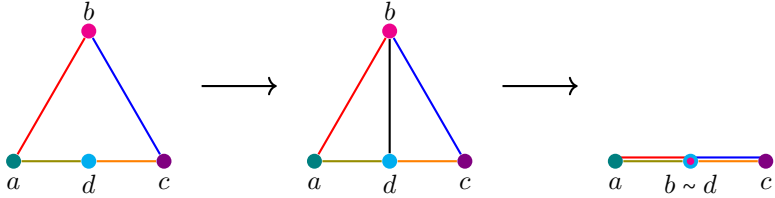
\begin{figure}[hbt!]
\centering
\begin{tikzpicture}
    \node (A) at (0,0) [circle, fill, inner sep=2pt, color=teal]{};
    \node at (0,-.3) {$a$};
    \node (C) at (2,0) [circle, fill, inner sep=2pt, color=violet]{};
    \node at (2,-.3) {$c$};
    \node (B) at (1,1.732) [circle, fill, inner sep=2pt, color=magenta]{};
    \node at (1,2) {$b$};
    \node (D) at (1,0) [circle, fill, inner sep=2pt, color=cyan]{};
    \node at (1,-.3) {$d$};
    \draw[thick, red] (A) -- (B);
    \draw[thick, blue] (B) -- (C);
    \draw[thick, orange] (C) -- (D);
    \draw[thick, olive] (D) -- (A);
    
   \draw[->, thick] (2.5,1) -- (3.5,1); 
    

\node (A) at (4,0) [circle, fill, inner sep=2pt, color=teal]{};
\node at (4,-.3) {$a$};
    \node (C) at (6,0) [circle, fill, inner sep=2pt, color=violet]{};
    \node at (6,-.3) {$c$};
    \node (B) at (5,1.732) [circle, fill, inner sep=2pt, color=magenta]{};
    \node at (5,2) {$b$};
    \node (D) at (5,0) [circle, fill, inner sep=2pt, color=cyan]{};
     \node at (5,-.3) {$d$};
    \draw[thick, red] (A) -- (B);
    \draw[thick, blue] (B) -- (C);
    \draw[thick, orange] (C) -- (D);
    \draw[thick, olive] (D) -- (A);  
    \draw[thick] (B) -- (D);


\draw[->, thick] (6.5,1) -- (7.5,1);


   \draw[thick, red]  (8,.05) -- (9, .05);
    \draw[thick, blue] (9, .05) -- (10, .05);
    \node (K) at (8,0) [circle, fill, inner sep=2pt, color=teal]{};
    \node (L) at (10,0) [circle, fill, inner sep=2pt, color=violet]{};
    \node (M) at (9,0) [circle, fill, inner sep=2pt, color=cyan]{};
    \node at (9,0) [circle, fill, inner sep=1pt, color=magenta]{};
 \draw[thick, olive] (K) -- (M);
    \draw[thick, orange] (M) -- (L);
    \node at (8, -.3) {$a$};
    \node at (9, -.3) {$b \sim d$};
    \node at (10,-.3) {$c$};
\end{tikzpicture}
\caption{Simplicial contractibility factorization.}\label{figure: thr contraction}
\end{figure}

We include the boundary of the subdivided Real 2-simplex into $\bar{\Delta}^2$ and then collapse through the 2-cells down to edge $adc$. Folding the subdivided 1-simplex to the left and gluing $a$ to $c$ produces a copy of $\simpotwo$ where the unit is included via $\eta_L$. If instead we fold the subdivided 1-simplex to the right and glue $c$ to $a$ we produce $\simpotwo$ where the unit has been included via $\eta_R$. Since the 1-simplex is contractible, both of these composites are null homotopic. We recover the wedge $\simpotwo \wedge_{D_2} \simpotwo$ at the center of the Hopf algebroid compatibility diagram in  \ref{diagram: thr compatibility} by gluing $a$ to $c$ in $\partial \bar{\Delta}^2$. Therefore we have that the maps $\mu_L \circ \delta '$ and  $\mu_R \circ \delta'$ are null homotopic since they factor through the contractible 1-simplex and the verification that the diagram in \ref{diagram: thr compatibility} is $D_2$-commutative up to homotopy is complete. 

Applying the functor $|A \otimes_{D_2} (-)|$ to the diagram in \ref{diagram: thr compatibility} yields the following diagram in $D_2$-spectra:

\begin{equation}
\begin{tikzcd}
A \arrow[d, "\eta_R"] & d\thr(A) \arrow[l, "\varepsilon \circ \pi"']  \arrow[r, "\varepsilon \circ \pi"] \arrow[d, "\delta'"] & A \arrow[d, "\eta_L"] \\
\thr(A) & \thr(A) \smsh_A \thr(A) \arrow[l, "\mu_R"'] \arrow[r, "\mu_L"] & \thr(A) \\
\thr(A) \smsh_{\mathbb{S}} \thr(A) \arrow[u, "\varphi"'] & \thr(A) \smsh_{\mathbb{S}} \thr(A) \arrow[u] \arrow[l, "\chi \smsh id"'] \arrow[r, "id \smsh \chi"] & \thr(A) \smsh_{\mathbb{S}} \thr(A). \arrow[u, "\varphi"]
\end{tikzcd}
\end{equation}

We note that $|A \otimes_{D_2} (\simpotwo \sqcup \simpotwo)|$ is the smash product of two copies of $\thr(A)$ as algebras over the sphere spectrum. Because the simplicial diagram was $D_2$-commutative up to homotopy, so too is the diagram in spectra and the proof that $\thr(A)$ is a Hopf algebroid in the $D_2$-homotopy category when $A$ is commutative is complete. 

\end{proof}

In the case of topological Hochschild homology, Angeltveit-Rognes show that the Hopf algebra structure on $\thh$ lifts to a Hopf algebra structure on the \bok\ spectral sequence. Specifically, Theorem 4.5 in \cite{angeltveit2005hopf} motivated the work undertaken in this section. A key step in the proof of their theorem involves using the simplicial Hopf algebra structure maps on $\thh$ to define maps on the spectral sequence. Such an approach is made possible because the \bok\ spectral sequence arises from a simplicial filtration of $\thh$. For this reason, we were careful to construct all of the Hopf algebroid structure maps simplicially in this section. Although we do not lift this structure to the Real \bok\ spectral sequence in this paper, we will return to this in future work. 

\begin{rem}
Van Niel \cite{danikathesis} classifies algebraic structures for the theory of $C_n$-twisted topological Hochschild homology discussed in Section \ref{section: twisted thh}. For $p$ prime and at least 5, the author proves that $\thh_{C_p}$ is a non-counital bialgebra in the equivariant stable homotopy category and that the twisted \bok\ spectral sequence of Theorem \ref{twistedbokss} inherits this structure. 
\end{rem}
  
\printbibliography
\end{document}